\documentclass[11pt]{article}
\usepackage{epsf,epsfig,color}
\usepackage{amssymb,latexsym, amsmath}
\usepackage{graphicx}%
\usepackage{soul}
\usepackage{pinlabel}

\newtheorem{lem}{Lemma}[section]
\newtheorem{theo}{Theorem}[section]
\newtheorem{pro}{Proposition}[section]
\newtheorem{cor}{Corollary}[section]
\newtheorem{con}{Conjecture}[section]

\newtheorem{prob}{Problem}[section]
\newtheorem{exa}{Example}[section]

\renewcommand{\theenumi}{\rm (\roman{enumi})}

\newcommand{\proof}
{{\noindent {\em Proof}.\quad}\setcounter{countclaim}{0}
\setcounter{countcase}{0}}
\newcommand{\proofend}{{\hfill$\Box$}}

\newcounter{countcase}

\newcounter{countclaim}
\def\inclaim{\addtocounter{countclaim}{1}
{\noindent {\bf Claim \thecountclaim}: }}

\newcounter{countfig}

\newcommand{\beeq}{\begin{equation}}
\newcommand{\eneq}{\end{equation}}

\newcommand{\beeqn}{\begin{eqnarray*}}
\newcommand{\eneqn}{\end{eqnarray*}}

\def \P {{\cal P}}
\def \OP {{\cal OP}}

\def \W {{\cal W}} 
\def \Z {\mathbb Z}
\def \Re {{\cal R}e}

\def \AO {{\cal AO}}
\def \repa {r}

\def \Gn {{\mathbb W}}

\def \iff {if and only if }

\newcommand {\relabel}[1] {\label{#1} \red{[*: #1]}}
\newcommand {\resection}[1] {\section{#1} }\newcommand {\rebibitem}[1] {\bibitem{#1} \red{[*: #1]}}
\def\relabel {\label} \def\resection {\section}\def\rebibitem {\bibitem}  

\begin{document}

\renewcommand{\theequation}{\thesection.\arabic{equation}}

\renewcommand{\theenumi}{\rm (\roman{enumi})}
\renewcommand{\labelenumi}{\rm(\roman{enumi})}

\baselineskip 0.6 cm

\title {
New expressions for order polynomials and chromatic polynomials
}

\author{
Fengming Dong\thanks{Email: fengming.dong@nie.edu.sg.} \\
\small Mathematics and Mathematics Education\\
\small National Institute of Education\\
\small Nanyang Technological University, 
Singapore 637616
}

\maketitle

\begin{abstract}
Let $G=(V,E)$ be a simple graph with 
$V=\{1,2,\cdots,n\}$
and $\chi(G,x)$ be its chromatic polynomial.
For an ordering $\pi=(v_1,v_2,\cdots,v_n)$ of elements of $V$,
let $\delta_G(\pi)$ be the number of $i$'s, where $1\le i\le n-1$,
with either $v_i<v_{i+1}$ or $v_iv_{i+1}\in E$.
Let $\W(G)$ be the set of
subsets $\{a,b,c\}$  
of $V$,
where $a<b<c$, which induces a subgraph 
with $ac$ as its only edge.
We show that $\W(G)=\emptyset$ \iff 
$(-1)^n\chi(G,-x)=\sum_{\pi} 
{x+\delta_G(\pi)\choose n}$, 
where the sum runs over all $n!$ orderings $\pi$ of $V$.
To prove this result, we establish an analogous 
result on order polynomials of posets  
and apply Stanley's work on 
the relation between chromatic polynomials
and  order polynomials. 
\end{abstract}


\resection {Introduction\relabel{intro}}

\subsection{Chromatic polynomials
\relabel{sec1-1}}

For a simple graph $G=(V,E)$, 
the {\it chromatic polynomial} of $G$  
is defined to be the polynomial  
$\chi(G,x)$ such that 
$\chi(G,k)$ counts 
the number of proper $k$-colourings of $G$ 
for any positive integer $k$ (for example, see~\cite{bir1912,bir1946,Dong2005,Read1968,RT1988,Whitney1932}).
This concept was first introduced by 
Birkhoff~\cite{bir1912} in 1912 in the hope of proving 
the four-color theorem
(i.e., $\chi(G,4)>0$ holds for any loopless planar graph $G$).
The study of chromatic polynomials is one of the
most active areas in graph theory
and many celebrated results on this topic
have been obtained (for example, see
\cite{bir1946,Dong2000,Jackson1993, Jackson2015, Royle2009, 
Sokal2001, Sokal2004, Stanley1973,  stanley1970, Thomassen1997, Thomassen2000,woo}).

One of the main purposes of this paper  
is to prove 
a new identity for $\chi(G,x)$ 
when $G$ satisfies a certain  
condition.
Assume that $V=[n]$, where $[n]=\{1,2,\cdots,n\}$.
For $u,v\in V$, define 
\begin{equation}\relabel{def-delta-G}
\delta_G(u,v)=
\left \{
\begin{array}{ll}
1, \qquad &u<v \mbox{ or } uv\in E;\\
0,&\mbox{otherwise}.
\end{array}
\right.
\end{equation}
Let $\P(V)$ denote the set of orderings of elements of $V$.
Obviously,  $|\P(V)|=n!$.
Define 
\begin{equation}\relabel{def-delta-11}
\Psi(G,x)=\sum_{\pi\in \P(V)}{x+\delta_G(\pi)\choose n},
\end{equation}
where for any $\pi=(u_1,u_2,\cdots,u_n)\in \P(V)$, 
\begin{equation}\relabel{def-delta-10}
\delta_G(\pi)=\sum_{1\le i\le n-1}\delta_G(u_i,u_{i+1}).
\end{equation}

Clearly the function $\Psi(G,x)$ depends on the structure 
of $G$ and also on the labeling of its vertices. 
For a bijection $\omega:V\rightarrow [n]$,
let $G_{\omega}$ denote the graph obtained from $G$ 
by relabeling each vertex $v$ in $G$ by $\omega(v)$.
Thus $G_{\omega}\cong G$
but it may be not true that 
$\Psi(G_{\omega},x)=\Psi(G,x)$.  
Hence, in this article,  
isomorphic graphs 
with different vertex labellings are 
considered to be different.  

For a graph $G=(V,E)$,  
where $V=[n]$, let 
$\W(G)$ be the set of $3$-element subsets 
$\{a,b,c\}$ of $V$ with $a<b<c$ 
such that $ac$ is the only edge in the subgraph of $G$ 
induced by $\{a,b,c\}$. 
Note that $\W(G)$ may be different from $\W(G_{\omega})$
for a bijection $\omega: V\rightarrow [n]$.

In Section~\ref{sect-chro}, we will prove the following 
result on $\chi(G,x)$.

\begin{theo}\relabel{cor-new1}
Let $G=(V,E)$ be a simple graph with 
$V=[n]$.  
Then 
\begin{equation}\relabel{cor-new1-eq1}
(-1)^n\chi(G,-x)
=\Psi(G,x)=\sum_{\pi\in \P(V)}{x+\delta_G(\pi)\choose n}
\end{equation}
\iff $\W(G)=\emptyset$.   
\end{theo}

To prove Theorem~\ref{cor-new1}, we will first establish
an analogous result on the order polynomial of 
$\bar D$ (i.e., Theorem~\ref{th-new2}), 
where $D$ is an acyclic digraph and 
$\bar D$ is the poset 
which is the reflexive transitive closure of $D$,
and apply Stanley's work on the relation 
between chromatic polynomials and order polynomials.

\subsection
{Order polynomials and strict 
order polynomials\relabel{sec1-2}}

In 1970, Stanley~\cite{stanley1970}
introduced the order polynomial and the strict 
order polynomial of a poset (i.e. partially ordered set).
Let $P$ be a poset on $n$ elements with a binary relation $\preceq $.
For $u,v\in P$, let $u\prec v$   
mean that 
$u\preceq v$ but $u\ne v$.
A mapping $\sigma: P\rightarrow [m]$ 
is said to be {\it order-preserving}   
(resp.,  {\it strictly order-preserving}) 
if $u\preceq v$ implies that $\sigma(u)\le \sigma(v)$
(resp.,  $u\prec v$ implies that $\sigma(u)<\sigma(v)$).
Let $\Omega(P,x)$ (resp.,  $\bar \Omega(P,x)$) 
be the function which 
counts the number of order-preserving (resp.,  strictly order-preserving) 
mappings $\sigma:P\rightarrow [m]$ 
whenever $x=m$ is a positive integer.
Both $\Omega(P,x)$ and $\bar \Omega(P,x)$
are polynomials in $x$ of degree $n$ 
(see Theorem 1 in \cite{stanley1970})
and are respectively called the {\it order polynomial} and 
the {\it strict order polynomial} of $P$.

An ordering $\pi=(v_1,v_2,\cdots,v_n)$ 
of the elements of $P$ is said to be {\it $P$-respecting}
if $v_i\prec v_j$ always implies that 
$i<j$ (i.e., $v_i$ appears before $v_j$ in $\pi$).
Let $\OP(P)$ be the set of $P$-respecting orderings $\pi$
of the elements of $P$.

Let $\omega$ be a fixed surjective order-preserving mapping
$\omega:P\rightarrow [n]$.
For a $P$-respecting ordering  
$\pi=(v_1,v_2,\cdots,v_n)$,
a ``decent" (resp. ``accent")  
means $\omega(v_i)>\omega(v_{i+1})$ 
(resp. $\omega(v_i)<\omega(v_{i+1})$)
for some $i$ with $1\le i\le n-1$.
Let $\kappa_{P}(\pi)$ (resp., $\bar \kappa_{P}(\pi)$)
denote the number of 
times when a ``decent" (resp. an ``accent")  
occurs in $\pi$.  
Clearly, $0\le \bar \kappa_P(\pi),\kappa_P(\pi)\le n-1$
and $\bar \kappa_P(\pi)+\kappa_P(\pi)=n-1$
for each $\pi\in \OP(P)$.
For an integer $s$  
with $0\le s\le n-1$,
let $w_s(P)$ (resp., $\bar w_s(P)$) 
be the number of $\pi\in \OP(P)$ 
with $\kappa_P(\pi)=s$
(resp., $\bar \kappa_P(\pi)=s$).

Stanley's Theorem 2 in \cite{stanley1970}  
gives the following 
interpretations for $\Omega(P,m)$ and $\bar \Omega(P,m)$.

\begin{theo}[Stanley~\cite{stanley1970}]\relabel{stan-th1}
For any integer $m\ge 1$,
\begin{align}
&\Omega(P,m)=\sum_{s=0}^{n-1} 
w_{s}(P){m+n-1-s\choose n}\mbox{ and } 
\bar \Omega(P,m)=\sum_{s=0}^{n-1} 
\bar w_{s}(P){m+n-1-s\choose n}.
\relabel{stanley-eq2}
\end{align}
\end{theo}
As $\kappa_P(\pi)+\bar \kappa_P(\pi)=n-1$ for each 
$\pi\in \OP(P)$, by applying Theorem~\ref{stan-th1}, 
it is not difficult to deduce that 
\begin{align}
\Omega(P,m)=\sum_{\pi\in \OP(P)} 
{m+\bar \kappa_P(\pi)\choose n}.
\relabel{stanley-eq3}
\end{align}
By Theorem~\ref{stan-th1}, 
a relation between  $\Omega(P,m)$ 
and $\bar \Omega(P,m)$ can also be deduced easily
and it appeared in Stanley's Theorem 3 in \cite{stanley1970}:
for any $m\in \Z^+$,
\begin{align}
\bar \Omega(P,m)=(-1)^n \Omega(P,-m).
\relabel{stanley-eq4}
\end{align}
From now on we focus on the order polynomial of 
a poset that is reflexive transitive closure 
of an acyclic digraph.  

A digraph $D=(V,A)$ is called {\it acyclic} 
if it does not contain any directed cycle.
Let $D$ be an acyclic digraph with $|V|=n$.
For convenience of notation, we simply  
assume that $V=[n]$. 
An ordering  
$\pi=(u_1,u_2,\cdots,u_n)$ of elements of $V$ 
is said to be {\it $D$-respecting}  
if $(u_i,u_j)\in A$ implies that $i<j$ holds  
(i.e., $u_i$ appears before $u_j$ in $\pi$).
Let $\OP(D)$ be the set of $D$-respecting 
orderings of elements of $V$. 
For example, for the digraphs in Figure~\ref{p1},
$\OP(D_i)$ has exactly three members 
given in Table~\ref{t1} for 
$i=1,2,3$.   

\begin{figure}[htbp]
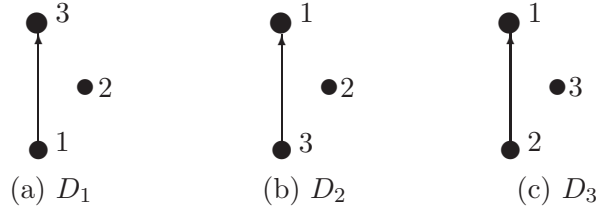

  \centering
\input p1.pic

(a) $D_1$ \hspace{2 cm} 
(b) $D_2$ \hspace{2 cm} 
(c) $D_3$

  \caption{Isomorphic digraphs $D_1,D_2$ and $D_3$}
\relabel{p1}
\end{figure}

\begin{table}[http]
 \centering
 
\begin{tabular}{c|c|c|c|c|c}
$\OP(D_1)$ & $\delta_{D_1}(\pi_i)$ 
&$\OP(D_2)$ & $\delta_{D_2}(\pi'_i)$  
&$\OP(D_3)$ & $\delta_{D_3}(\pi''_i)$ \\ \hline
$\pi_1=(2,1,3)$  & $1$  & $\pi'_1=(2,3,1)$ & $2$
& $\pi''_1=(3,2,1)$ & $1$\\ \hline
$\pi_2=(1,2,3)$  & $2$  & $\pi'_2=(3,2,1)$ & $0$
& $\pi''_1=(2,3,1)$ & $1$\\ \hline
$\pi_3=(1,3,2)$  & $1$  & $\pi'_3=(3,1,2)$ & $2$
& $\pi''_1=(2,1,3)$ & $2$\\
\end{tabular}

 \caption{Members of $\OP(D_i)$ and values
$\delta_{D_i}(\pi)$ for $\pi\in \OP(D_i)$}
 \relabel{t1}
\end{table}

Clearly, an ordering $\pi$ of elements of $V$
is $D$-respecting \iff it is $\bar D$-respecting.
Thus $\OP(D)=\OP(\bar D)$.

For $a,b\in \Z^+$,  
let $\bar \kappa(a,b)=1$ if $a<b$, and 
$\bar \kappa(a,b)=0$ otherwise. 
For an ordering $\pi=(a_1,a_1,\cdots,a_n)$ 
of $n$ different numbers in $\Z^+$, let 
$$
\bar \kappa(\pi)=\sum_{i=1}^{n-1}\bar\kappa(a_i,a_{i+1}).
$$
Thus $\bar \kappa(\pi)$ is actually the number of times 
when an ``accent" occurs  
in the ordering $\pi$.
Note that the definition of $\bar \kappa(\pi)$ 
is only related to the numbers in the ordering $\pi$
and has no relation with $D$.  

Let $\Re(D)=\{(a,b)\in A: a>b\}$. 
Assume that $\Re(D)=\emptyset$.
As $V=[n]$, 
this assumption is equivalent to 
a surjective mapping $\omega:V\rightarrow [n]$ 
with the property that 
$(u,v)\in A$ 
implies $\omega(u)<\omega(v)$.
Observe that for any $\pi\in \OP(D)$,
$\bar \kappa(\pi)=\bar \kappa_{\bar D}(\pi)$ holds.
Thus, by (\ref{stanley-eq3}),  
$\Omega(\bar D,m)$ has the following 
expression in terms of $\bar \kappa(\pi)$
under the assumption that $\Re(D)=\emptyset$:
\begin{align}
\Omega(\bar D,m)
=\sum_{\pi\in \OP(D)}
{m+\bar \kappa(\pi)\choose n}.
\relabel{eq7}
\end{align}
Note that 
if $\Re(D)\ne \emptyset$, 
(\ref{eq7}) may be not true,  
unless $\bar \kappa(\pi)$ is 
replaced by another suitable function.  
In the following, we remove the assumption that $\Re(D)=\emptyset$
and replace $\bar \kappa(\pi)$ by a new
function $\delta_D(\pi)$. 
We will see for which labellings of vertices of $D$ 
an identity analogous to (\ref{eq7}) holds
even if $\Re(D)\ne \emptyset$.  

\subsection
{A new function $\Psi(D,x)$ for an acyclic digraph $D$\relabel{sec1-3}}

Let $D=(V,A)$ be an acyclic digraph with $V=[n]$. 
For $a,b\in V$, define 
\begin{equation}\relabel{def-delta-0}
\delta_D(a,b)=
\left \{
\begin{array}{ll}
1, \qquad &\mbox{either }a<b \mbox{ or } (a,b)\in A;\\
0,&\mbox{otherwise}.
\end{array}
\right.
\end{equation}
Clearly $\kappa(a,b)\le \delta_D(a,b)$
for every pair of members $a$ and $b$ of $V$.
When $\Re(D)=\emptyset$, 
$(a,b)\in A$ implies that $a<b$.
Thus, in this case,  
$\delta_D(a,b)=\bar \kappa(a,b)$ holds
for every pair of numbers $a$ and $b$ in $V$,
no matter whether $(a,b)\in A$ or not.
However, when $\Re(D)\ne \emptyset$, 
for each $(a,b)\in A$ with $a>b$,
we have $\delta_D(a,b)=1$ and $\bar \kappa(a,b)=0$.

Let $\Psi(D,x)$ be the function defined below:
\begin{equation}\relabel{stan-eq2}
\Psi(D,x)=\sum_{\pi\in \OP(D)}{x+\delta_D(\pi)\choose n},
\end{equation}
where for any $\pi=(a_1,a_2,\cdots,a_n)\in \OP(D)$, 
\begin{equation}\relabel{def-gamma}
\delta_D(\pi)=\sum_{i=1}^{n-1} \delta_D(a_i,a_{i+1}).
\end{equation}
Note that $\Psi(D,x)$ is a function defined on an
acyclic digraph $D=(V,A)$ 
with $V$ a linearly ordered set of $n$ vertices  
and its definition does not rely on 
a fixed mapping $\omega:V\rightarrow [n]$
with the property that $(v_i,v_j)\in A$ 
implies  $\omega(v_i)<\omega(v_j)$.

Clearly, if $\Re(D)=\emptyset$, 
then $\delta_D(\pi)=\bar \kappa(\pi)$ holds 
for every $\pi\in \OP(D)$,
and thus  (\ref{eq7}) and (\ref{stan-eq2})
imply the following conclusion.

\begin{pro}
\relabel{stan1}
Let $D=([n],A)$ be an acyclic digraph.   
If $\Re(D)=\emptyset$, then 
\begin{equation}\relabel{stan1-eq1}
\Omega(\bar D,x)=\Psi(D,x)=
\sum_{\pi\in \OP(D)}{x+\delta_D(\pi)\choose n}.
\end{equation}
\end{pro}

If  $\Re(D)\ne \emptyset$, it is possible 
that $\delta_D(\pi)\ne \bar \kappa(\pi)$ 
for some $\pi\in \OP(D)$,
and thus it is possible that 
$\Omega(\bar D,x)\ne \Psi(D,x)$.
For example, for the isomorphic digraphs 
$D_1,D_2$ and $D_3$ in Figure~\ref{p1}, 
by the data in Table 1,
we have 
\begin{equation}\relabel{stan1-eq3}
\Psi(D_1,x)=\Psi(D_3,x)={x+2\choose 3}+2{x+1\choose 3}
\ne \Psi(D_2,x)=2{x+2\choose 3}+{x\choose 3}.
\end{equation}
As $\Re(D_1)=\emptyset$, by Proposition~\ref{stan1},
we have $\Psi(D_3,x)=\Psi(D_1,x)=\Omega({\bar D_1},x)
=\Omega({\bar D_3},x)$.
But $\Psi(D_2,x)\ne \Psi(D_1,x)=\Omega({\bar D_1},x)=\Omega({\bar D_2},x)$.

Notice that $\Re(D_3)\ne \emptyset$,
although $\Psi(D_3,x)=\Omega({\bar D_3},x)$.
Thus, 
$\Psi(D,x)=\Omega(\bar D,x)$ does not imply 
$\Re(D)=\emptyset$.
The main aim of this article is to determine exactly when
the identity $\Omega(\bar D,x)=\Psi(D,x)$ holds.   

Let $D=(V,A)$ be an acyclic digraph, where $V=[n]$. 
For distinct $a,b\in V$, write $a\prec_D b$ if 
there exists a directed path in $D$ connecting 
from $a$ to $b$,
and $a\not \prec_D b$ otherwise. 
Write $a\not \approx_D b$ if $a\not \prec_D b$ 
and $b\not \prec_D a$.
Let $\W(D)$ be the set of 3-element subsets $\{a,b,c\}$ 
of $V$ with $a<b<c$ such that 
$(c,a)\in A$ but $a\not \approx_D b$ and $c\not \approx_D b$.
Observe that if  $(c,a)\in A$, 
then $b\prec_D c$ implies that $b\prec_D a$,
and $a\prec_D b$ implies that $c\prec_D b$.
Thus, for $\{a,b,c\}\subseteq V$ with $a<b<c$ and $(c,a)\in A$, 
$\{a,b,c\}\in \W(D)$ 
\iff $c\not\prec_D b$ and $b\not\prec_D a$.

For example, for the digraphs $D_1, D_2$ and $D_3$ in Figure~\ref{p1},
only $\W(D_2)$ is not empty, and 
for the digraph $D$ in Figure~\ref{f2} 
on Page~\pageref{f2},  
$\W(D)$ has exactly one 
member $\{2,3,5\}$.

Clearly, $\Re(D)=\emptyset$ implies that 
$\W(D)=\emptyset$. 
But the converse does not hold.
In Section~\ref{identity}, 
we will show that if $\W(D)=\emptyset$,
then there exists $D'$ obtained from $D$ by relabeling 
vertices in $D$ such that $\Re(D')=\emptyset$
and $\Psi(D,x)=\Psi(D',x)$.
By Proposition~\ref{stan1}, we have 
$\Psi(D',x)=\Omega({\bar D'},x)=\Omega(\bar D,x)$.
Thus we establish the following result.

\begin{theo}\relabel{th-new1}
Let $D = ([n],A)$ be an acyclic graph and 
$\W(D)$ be defined as above. 
If $\W(D) = \emptyset$, then  
$\Psi(D,x)=\Omega(\bar D,x)$ holds.   
\end{theo}

The converse of Theorem~\ref{th-new1} also holds,
as stated in the following result.

\begin{theo}\relabel{th-new2} 
Let $D = ([n],E)$ be an acyclic graph, where $n \ge 3$. 
Then   
\begin{equation}\relabel{main-th1-eq2}
\Psi(D,x)-\Omega(\bar D,x)=\sum_{i=0}^{n-3}d_i{x+i\choose n-2},
\end{equation}
where $d_0,d_1,\cdots,d_{n-3}$ are non-negative integers.
Furthermore, 
$d_i=0$ for every $i=0,1,\cdots,n-3$ 
\iff $\W(D)=\emptyset$.
\end{theo}

Clearly, Theorem~\ref{th-new2} implies that 
$\Psi(D,x)=\Omega(\bar D,x)$ \iff $\W(D)=\emptyset$.
To prove Theorem~\ref{th-new2} in Section~\ref{sec3}, we will first compare $\Psi(D,x)$ with $\Psi(D_{a\rightarrow r},x)$,
where $D_{a\rightarrow r}$ is the digraph obtained from 
$D$ by relabeling vertex $a$ by by a suitable number $r$.
The new digraph $D_{a\rightarrow r}$ has the property that 
$\W(D_{a\rightarrow r})=\W(D)-\{W\in \W(D): a\in W\}$
and $\Psi(D,x)-\Psi(D_{a\rightarrow r},x)
=\sum\limits_{i=0}^{n-3} d_i{x+i\choose n-2}$, 
where $d_i\ge 0$ for all $i$, and 
$d_0+\cdots+d_{n-3}=0$ \iff $\W(D)=\emptyset$.

While Theorem~\ref{th-new1} is implied by 
Theorem~\ref{th-new2}, 
the derivation of Theorem~\ref{th-new2}  
is independent of Theorem~\ref{th-new1}.  
For a special case, the numbers $d_i$ in Theorem~\ref{th-new2} are given an interpretation
(see Proposition~\ref{pro5-5}).

Let $\AO(G)$ be the set of acyclic orientations of $G$.
The expression (1) in \cite{Stanley1973}  
gives a relation between 
$\chi(G,x)$ and 
$\bar \Omega(\bar D,x)$:
\begin{align}
\chi(G,x)=\sum_{D\in \AO(G)} \bar \Omega(\bar D,x).
\relabel{stanley-eq5}
\end{align}
 
Thus, (\ref{stanley-eq3}), (\ref{stanley-eq4}) 
and (\ref{stanley-eq5}) imply the following result.

\begin{theo}[Stanley \cite{Stanley1973}] 
\relabel{stan2}
Let $G=(V,E)$ be a simple graph. 
Then 
\begin{equation}\relabel{stan2-eq2}
(-1)^{|V|}\chi(G,-x)
=\sum_{D\in \AO(G)} \Omega(\bar D,x).
\end{equation}
\end{theo}

Note that for each $D\in \AO(G)$, 
determining $\Omega(\bar D,x)$ by (\ref{eq7})
is based on a relabeling of vertices 
such that $a<b$ holds
for each arc $(a,b)$ in $D$.
Thus, the summation of (\ref{stan2-eq2}) cannot be 
replaced by a summation over all $|V|!$ orderings
of elements of $V$ 
if the labeling of elements of $V$ is fixed,
although the union of $\OP(D)$'s for all $D\in \AO(G)$ 
is exactly the set of all $|V|!$ orderings
of elements of $V$. 
This is another motivation for 
extending (\ref{eq7}) to an analogous expression  
with an arbitrary relabeling
of vertices in $D$ and the result can be applied 
to express $\chi(G,x)$ 
as the summation over all 
$|V|!$ orderings of elements of $V$.

Applying Theorems~\ref{th-new2} and~\ref{stan2},
we can prove Theorem~\ref{cor-new1} in Section~\ref{sect-chro}.

\section{Proof of Theorem~\ref{th-new1}
\relabel{identity}}

Let $D=(V,A)$ be an acyclic digraph with vertex 
set $V$, where $V=[n]$.  
In this section, we shall show that 
$\Psi(D,x)=\Omega(\bar D,x)$ whenever $\W(D)=\emptyset$.

For $S\subseteq V$, let $D[S]$ be the subdigraph 
of $D$ induced by $S$. 
For $u\in V$,  
$u$ is called a {\it sink} of $D$
if $F_D(u)=\emptyset$, 
where $F_D(u)=\{v: (u,v)\in A\}$.
We first define a bijection $L:V\rightarrow [n]$ 
by the following algorithm:

\noindent {\bf Algorithm A}:
\begin{enumerate}
\item[Step 1.] Set $S:=V$;
\item[Step 2.] Let $u$ be the largest number 
among all sinks of $D[S]$; 
\item[Step 3.] Set $L(u):=|S|$ and $S:=S\setminus \{u\}$;
\item[Step 4.] If $S\ne \emptyset$, go to Step 2; 
otherwise, output $L(v)$ for all $v\in V$.
\end{enumerate}
The bijection $L$ defined above will be written as $L_D$ when
there is a possibility of confusion.  

\begin{exa}\relabel{ex1}
If $D$ is the acyclic digraph in Figure~\ref{f2}, then 
$$
L(3)=5, L(2)=4,L(5)=3, L(4)=2,L(1)=1.
$$
\end{exa}

\begin{figure}[h!]
\centering 
\scalebox{0.9}{\input{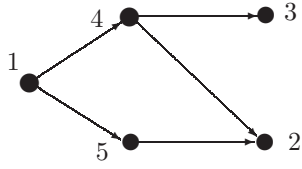}}
\caption{An acyclic digraph}
\relabel{f2}
\end{figure}

Recall that for distinct $u,v\in V$, 
$u\prec_D v$ if $D$ has a directed path from $u$ to $v$;
and for  $u\in V$, 
$R_D(u)$ (or simply $R(u)$) 
denote the set $\{v\in V: u\prec_D v\}$.
Let $R_D[u]=\{u\}\cup R_D(u)$.
Then $u\in R_D[u]$ but $u\notin R_D(u)$.

By definitions of $\prec_D$ and $L_D$, 
we have the following basic properties of  
$\prec_D$ and $L_D$.

\begin{pro}\relabel{pro2-0} 
Let $a,b$ and $c$ be distinct vertices in $D$. 
\begin{enumerate}
\renewcommand{\theenumii} {\mbox{\roman{enumii}}}
\item\relabel{pro2-1-no1} 
If $a\prec_D b$ and $b\prec_D c$, then $a\prec_D c$.

\item\relabel{pro2-1-no2} If $a\prec_D b$, then 
$L(a)<L(b)$.
\end{enumerate}
\end{pro}

For distinct vertices $b,c$ in $D$, let 
$N_D[c,b]=\{c'\in R_D[c]\setminus R_D[b]:
\forall y\in R_D(c)\cap R_D(b), L(c')<L(y)\}$.  

\begin{exa}\relabel{ex2}
For the digraph $D$ in Figure~\ref{f2},
$N_D[5,3]=\{5,2\}$ and $N_D[5,4]=\{5\}$.
\end{exa}

\begin{pro}\relabel{pro2-01} 
Let $b$ and $c$ be distinct vertices in $D$
with $c\notin R_D(b)$.  Then  
\begin{enumerate}
\item  $c\in N_D[c,b]$;
\item when $R_D(c)\subseteq R_D(b)$, 
 $N_D[c,b]=\{c\}$ holds.
\end{enumerate}
\end{pro}

\proof (i). Clearly $c\in R_D[c]\setminus R_D[b]$.
As $R_D(c)\cap R_D(b)\subseteq R_D(c)$, 
we have $L(c)< L(y)$ for all 
$y\in R_D(c)\cap R_D(b)$
by Proposition~\ref{pro2-0}~\ref{pro2-1-no2},  
implying that $c\in N_D[c,b]$.
Thus (i) holds.

(ii). By the result in (i), $c\in N_D[c,b]$.
As $R_D(c)\subseteq R_D(b)$, $R_D[c]\setminus R_D[b]=\{c\}$.
Thus (ii) holds.
\proofend

For an non-empty finite set $S$ of $\Z^+$, 
let $\min S$ and $\max S$ denote the minimum value 
and the maximum value of $S$ respectively. 
In case of any confusion, 
 $\min S$ and $\max S$ are respectively written as 
 $\min(S)$ and $\max(S)$.

The bijection $L_D:V\rightarrow \{1,2,\cdots,n\}$ 
has the following  property.

\begin{pro}\relabel{pro2-1} 
Let $a,b$ and $c$ be distinct vertices in $D$. 
\begin{enumerate}
\renewcommand{\theenumii} {\mbox{\roman{enumii}}}

\item\relabel{pro2-1-no3}  
If $c\not\approx_D b$, then 
$L(c)<L(b)$ \iff 
$\min (N_D[c,b])< \min (N_D[b,c])$;   

\item\relabel{pro2-1-no4}   
If $c\not \approx _D b$,  $L(c)<L(b)$ and $b<c$,
then there exist $a,c'\in R_D[c]\setminus R_D[b]$  
such that $\{a,b,c'\}\in \W(D)$;
 
\item\relabel{pro2-1-no5}
If $\W(D)=\emptyset$, $b<c$ and $c\not \prec _D b$,
then  $L(b)<L(c)$.
\end{enumerate}
\end{pro}

\proof 
\ref{pro2-1-no3}. 
Assume that $c\not \approx_D b$.
It suffices to prove  
that if
$\min (N_D[c,b])< \min (N_D[b,c])$, 
then $L(c)<L(b)$, 
as exchanging $b$ and $c$ yields that 
if $\min (N_D[b,c])< \min (N_D[c,b])$, 
then $L(b)<L(c)$.  

By Proposition~\ref{pro2-01} (i),
$c\in N_D[c,b]$ and 
$b\in N_D[b,c]$.
Let $c_0=\min (N_D[c,b])$. 
By Proposition~\ref{pro2-0}~\ref{pro2-1-no2},
$L(c)\le L(c_0)$.

Let $S'$ be the set of sinks of $D$ and 
let $w=\max S'$.
Then $L(w)=|V|$. 
Now we want to prove the two following claims  
under the assumption that $c_0<\min (N_D[b,c])$.  

\noindent {\bf Claim 1}: $w\ne c_0$.

Assume that $w= c_0$.
As  $L(c_0)=|V|$, 
$c_0$ is the largest sink of $D$.
Note that $S'\cap R_D[b]\ne \emptyset$. 
Let $b_0=\max (S'\cap R_D[b])$.
As $c_0\in R_D[c]\setminus R_D[b]$, we have $b_0\ne c_0$
and so $b_0<c_0$ and $L(b_0)<L(c_0)=|V|$.
As $b_0<c_0<\min (N_D[b,c])$ and $b_0\in R_D[b]$, 
we have $b_0\in R_D[b]\setminus N_D[b,c]$.
By the assumption on $N_D[b,c]$, 
$b_0\in R_D[b]\setminus N_D[b,c]$ implies that 
$b_0\in R_D(c)\cap R_D(b)$ 
or $L(b_0)>L(y)$ for some $y\in R_D(c)\cap R_D(b)$.
Thus $L(b_0)\ge L(y)$ for some $y\in R_D(c)\cap R_D(b)$.
As $L(c_0)<L(y)$ 
for all $y\in R_D(c)\cap R_D(b)$, 
we have $L(c_0)<L(b_0)$, a contradiction.


\noindent {\bf Claim 2}:
$L(c_0)<L(b)$.

This claim is trivial when $|V|=2$. 
Now assume  $|V|\ge 3$ and that this claim fails.  
Thus $L(b)<L(c_0) \le |V|$.

By Claim 1, $w\ne c_0$.
Then $L(c)\le L(c_0)<L(w)=|V|$.  
As Claim 1 holds for $D-w$, by induction,  
and 
$$
\min (N_{D-w}[c,b])=\min (N_{D}[c,b])=c_0
<\min (N_{D}[b,c])=\min (N_{D-w}[b,c]),
$$
we have 
$L_{D-w}(c_0)<L_{D-w}(b)$.
Since 
$L_{D-w}(c_0)=L_D(c_0)$ and $L_D(b)=L_{D-w}(b)$,
we have $L_{D}(c_0)<L_{D}(b)$,
a contradiction. 
Thus Claim 2 holds.

As $L(c)\le L(c_0)$, 
Claim 2 implies $L(c)<L(b)$
under the condition that 
$\min(N_D(c,b))<\min(N_D(b,c))$.
Thus \ref{pro2-1-no3} holds.

\ref{pro2-1-no4}. 
Assume that 
$b\not \approx _D c$, $b<c$ and $L(c)<L(b)$.
By \ref{pro2-1-no3}, $\min (N_D[c,b])<\min (N_D[b,c])$.
Let $c_1=\min(N_D[c,b])$.
Then $c_1<\min (N[b,c])\le b<c$.
As $c_1\in N_D[c,b]\subseteq R_D[c]$, there is a path in $D$ 
from $c$ to $c_1$:
$c\rightarrow a_1\rightarrow \cdots \rightarrow a_k$,
where $a_k=c_1$ and $a_i\rightarrow a_{i+1}$ is  
short for $(a_i,a_{i+1})\in A$.  
As $a_k=c_1<b<c$, there exists $i: 1\le i\le k-1$ 
such that $a_i>b>a_{i+1}$. 
As $c_1\in N_D[c,b]\subseteq R_D[c]\setminus R_D[b]$, 
we have $a_i,a_{i+1}\in R_D[c]\setminus R_D[b]$,   
implying that $b\not\approx_D a_i$ and 
$b\not\approx_D a_{i+1}$.
Thus $\{a_{i+1},b,a_i\}\in \W(D)$ and the result holds. 

\ref{pro2-1-no5}.
Assume that $\W(D)=\emptyset$, $b<c$ and
$c\not \prec_D b$. 
If $b\prec _D c$, then Proposition~\ref{pro2-0}~\ref{pro2-1-no2}
implies that $L(b)<L(c)$.
Now assume that $b\not\prec _D c$.
Thus  $b\not \approx_D c$.
As $\W(D)=\emptyset$ and $b<c$,
by \ref{pro2-1-no4}, we have $L(b)<L(c)$
in this case. 
\proofend

Let $D_L$ be the digraph obtained from $D$ 
by relabeling each vertex $y$ in $D$ as $L(y)$.
Clearly, $D_L$ is isomorphic to $D$
and Proposition~\ref{pro2-0}~\ref{pro2-1-no2}
implies that $\Re(D_L)=\emptyset$.
By Proposition~\ref{stan1}, $\Psi(D_L,x)=\Omega({\bar D_L},x)=\Omega(\bar D,x)$.

For  $\pi=(a_1,a_2,\cdots,a_n)\in \OP(D)$,
let $L(\pi)=(L(a_1),L(a_2),\cdots,L(a_n))$.

\begin{pro}\relabel{pro2-2}  
Let $\pi=(a_1,a_2,\cdots,a_n)\in \OP(D)$.    
If $\W(D)=\emptyset$,  then 
\begin{enumerate}
\item\relabel{pro-2-no1}   
$\delta_{D}(a_i,a_{i+1})=\delta_{D_L}(L(a_i),L(a_{i+1}))$
holds for $i=1,2,\cdots,n-1$; and

\item\relabel{pro-2-no2}
$\delta_D(\pi)=\delta_{D_L}(L(\pi))$ holds.
\end{enumerate}
\end{pro}

\proof \ref{pro-2-no1}.
As $\pi=(a_1,a_2,\cdots,a_n)\in \OP(D)$,
we have $a_{i+1}\not \prec a_i$. 
Thus, either $a_i\prec_D a_{i+1}$ or 
$a_i\not \approx _D a_{i+1}$.

First consider the case that $a_i\prec_D a_{i+1}$.
As $\pi=(a_1,a_2,\cdots,a_n)\in \OP(D)$, 
if $a_{j_1}\rightarrow a_{j_2}\rightarrow \cdots 
\rightarrow a_{j_k}$ is a path in $D$, 
then $j_1<j_2<\cdots<j_k$.
Thus $a_i\prec_D a_{i+1}$ implies that 
$(a_i,a_{i+1})\in A$, 
and so $\delta_D(a_i,a_{i+1})=1$.
As $(a_i,a_{i+1})\in A$, we have 
$(L(a_i),L(a_{i+1}))\in A(D_L)$
and so $\delta_{D_L}(L(a_i), L(a_{i+1}))=1$.

Now assume that 
$a_i\not \approx _D a_{i+1}$.
As $\W(D)=0$, by Proposition~\ref{pro2-1} \ref{pro2-1-no5},
if $a_i<a_{i+1}$ then $L(a_i)<L(a_{i+1})$;
if $a_{i+1}<a_{i}$ then $L(a_{i+1})<L(a_{i})$.    
As $a_i\not \approx _D a_{i+1}$, 
we have $(a_i,a_{i+1})\notin A(D)$
and $(L(a_i),L(a_{i+1}))\notin A(D_L)$.
By definition of $\delta_D(a_i,a_{i+1})$, 
$\delta_{D_L}(L(a_i),L(a_{i+1}))=\delta_D(a_i,a_{i+1})$
holds in this case.

Thus \ref{pro-2-no1} holds.
By the result in \ref{pro-2-no1}, 
\ref{pro-2-no2} follows
directly from the definition of $\delta_D(\pi)$.
\proofend

\begin{cor}\relabel{cor-2-1}
If $\W(D)=\emptyset$, then 
$\Psi(D,x)=\Psi(D_L,x)$
\end{cor}

\proof
Note that 
$\pi\in \OP(D)$ \iff $L(\pi)\in \OP(D_L)$.
Thus
$$
\OP(D_L)=\{L(\pi): \pi\in \OP(D)\}.
$$
By Proposition~\ref{pro2-2}~\ref{pro-2-no2}, 
$\delta_D(\pi)=\delta_{D_L}(L(\pi))$ holds 
for each $\pi\in \OP(D)$. 
By definition of $\Psi(D,x)$, 
$\Psi(D,x)=\Psi(D_L,x)$ holds.
\proofend

Since $\Re(D_L)=\emptyset$, 
Proposition~\ref{stan1} implies that 
$\Psi(D_L,x)=\Omega({\bar D_L},x)=\Omega(\bar D,x)$.
Thus Theorem~\ref{th-new1} follows from 
Corollary~\ref{cor-2-1}.

\section{Proof of Theorem~\ref{th-new2}
\relabel{sec3}    
}

In this section, we assume 
that $D=(V,A)$ is an acyclic digraph with $V\subset \Z^+$ 
and $|V|=n$, 
where $n\ge 3$.
For $a\in V$ and $r\in \Z^+\setminus V$, 
let $D_{a\rightarrow r}$ be the digraph obtained from $D$ 
by relabeling $a$ by $r$. 
We will compare  
$\Psi(D,x)$ with $\Psi(D_{a\rightarrow r},x)$
and apply the result on 
$\Psi(D,x)-\Psi(D_{a\rightarrow r},x)$
to prove Theorem~\ref{th-new2}.

Clearly, if $V=[n]$, then $r\ge n+1$ and 
the vertex set of $D_{a\rightarrow r}$ 
is $([n]\setminus \{a\})\cup \{r\}$ which is 
no longer $[n]$.
Thus, for the purpose of comparing 
$\Psi(D,x)$ with $\Psi(D_{a\rightarrow r},x)$,
in this section the vertex set $V$ is allowed 
to be any subset of $\Z^+$  
and it is possible that $V\ne [n]$.


Note that if $V=\{v_1,v_2,\cdots,v_n\}$ 
with $1\le v_1<v_2<\cdots<v_n$,
then $\Psi(D,x)=\Psi(D',x)$ holds, 
where $D'$ is obtained from $D$ by 
relabeling each $v_i$ by $i$.
So the function $\Psi(D,x)$ is not affected 
even if  $V\ne [n]$.

\subsection{
Relabel a vertex in $D$ by a sufficiently large number
\relabel{sec3-1}
}
  
Define 
\begin{equation}\relabel{main-th1-eq1}
\Delta(D,z)=\sum_{\pi\in \OP(D)}z^{\delta_D(\pi)}.
\end{equation}
By definitions of $\Psi(D,x)$ and $\Delta(D,z)$, 
for any two acyclic digraphs $D_1$ and $D_2$ 
of the same order, 
$\Delta(D_1,z)=\Delta(D_2,z)$ 
\iff $\Psi(D_1,x)=\Psi(D_2,x)$. 

In this subsection, we always assume that 
$a$ is a fixed vertex in $D$ and 
$m$ is a number in $\Z^+\setminus V$ with 
$m>y$ for all $y\in V\setminus R_D[a]$.
We compare $\Delta(D,z)$ with 
$\Delta(D_{a\rightarrow m},z)$
under this assumption. 
This result will be applied in the next subsection 
for relabeling vertex $a$ by a suitable number $\repa$ 
so that $D$ can be replaced by $D_{a\rightarrow \repa}$ 
for the purpose of proving Theorem~\ref{th-new2}.

\subsubsection{
A function $\Delta_{D,\pi_0}(z)$
\relabel{sec3-1-1}}

Let 
$\pi_0=(a_1,a_2,\cdots,a_{n-1})$ be a fixed member
of $\OP(D-a)$, where 
$D-a$ is the digraph obtained from $D$ by removing 
vertex $a$.
Let $\OP(D,\pi_0)$ be the set of those members $\pi\in \OP(D)$ 
such that $\pi-a=\pi_0$,
where $\pi-a$ is obtained from $\pi$ by removing $a$.
For example, if $\pi=(2,1,3,4)$, 
then $\pi-2=(1,3,4)$. 
Observe that $(a_1,a_2,\cdots,a_{n-1},a)\in \OP(D,\pi_0)$
\iff $a$ is a sink of $D$, 
and $(a_1,\cdots, a_i,a,a_{i+1},\cdots,a_{n-1})\in \OP(D,\pi_0)$
\iff $(a_j,a)\notin A$ for all $j=i+1,\cdots,n-1$
and $(a,a_j)\notin A$ for all $j=1,\cdots,i$.

A vertex $u$ of $D$ is called a {\it source} if 
$(v,u)\not\in A$ for all $v\in V$.
Throughout this section, let $s$ and $t$ be 
the two numbers 
defined below:
\begin{enumerate}
\item let $s=0$ if $a$ is a source of $D$, 
and let $s=\max \{1\le k\le n-1: (a_k,a)\in A\}$ otherwise;
\item 
let $t=n$ if $a$ is a sink of $D$,
and let $t=\min \{1\le k\le n-1: (a,a_k)\in A\}$ otherwise.
\end{enumerate}
If $s=0$ or $t=n$, then clearly $s<t$. 
Otherwise, 
$(a_s,a)\in A$ and $(a,a_t)\in A$ imply that 
$a_s\prec _D a_t$,
and so $s<t$
by the assumption that $\pi_0\in \OP(D-a)$.
Hence we always have $s<t$.

By definition of $\OP(D)$ and the assumptions on $s$ 
and $t$,    
we have 
\begin{equation}\relabel{OP-a}
\OP(D,\pi_0)=\{(\cdots,a_i,a,a_{i+1},\cdots):
s\le i\le t-1\}.
\end{equation}
For $\pi\in \OP(D)$, let $\pi_{a\rightarrow m}$ 
be the ordering  
obtained from $\pi$ by replacing $a$ by $m$.
Then, 
\begin{equation}\relabel{OP-m}
\OP(D_{a\rightarrow m},\pi_0)=\{\pi_{a\rightarrow m}: \pi\in \OP(D,\pi_0)\}
=\{(\cdots,a_i,m,a_{i+1},\cdots):
s\le i\le t-1\}.
\end{equation}
Define 
\begin{equation}\relabel{def-delta-1} 
\Delta_{D,\pi_0}(z)
=\sum_{\pi\in \OP(D,\pi_0)}z^{\delta_D(\pi)-\delta_{D-a}(\pi_0)}.
\end{equation}
By (\ref{OP-a}), we have  
\begin{equation}\relabel{def-delta-2} 
\Delta_{D,\pi_0}(z)
=\sum_{s\le i\le t-1}
z^{\delta_D(a_i,a)+\delta_D(a,a_{i+1}) 
-\delta_{D}(a_i,a_{i+1})},
\end{equation}
where the following numbers  are assumed in case
that $s=0$ or $t=n$:
\begin{equation}\relabel{def-delta-3} 
1=\delta_D(a_0,a_1)
=\delta_D(a_0,a)
=\delta_D(a_{n-1},a_n)
=\delta_D(a, a_{n}).
\end{equation}

\subsubsection{Expression for   
$\Delta(D,z)-\Delta(D_{a\rightarrow m},z)$ 
\relabel{sec3-1-2}
}

Let $U_1$ and $U_2$ be the two disjoint subsets of 
$\{i: s+1\le i\le t-2\}$ defined below:
\begin{equation}
\left \{
\begin{array}{l}
U_1=\{s+1\le i\le t-2: a_i>a>a_{i+1}\},\\
U_2=\{s+1\le i\le t-2: a_i<a<a_{i+1}\}.
\end{array}
\right.\relabel{le3-1-beq1}
\end{equation}

\begin{lem}\relabel{le3-1}
(i) $\Delta_{D,\pi_0}(z)$ 
has the following expression:
\begin{eqnarray}
\Delta_{D,\pi_0}(z)
&=&z^{1+\delta_D(a,a_{s+1})
-\delta_D(a_s,a_{s+1})}+z^{1+\delta_D(a_{t-1},a)
-\delta_D(a_{t-1},a_{t})}
+\sum_{i\in U_1}z^{-\delta_D(a_i,a_{i+1})}\nonumber \\
& &
+\sum_{i\in U_2}z^{2-\delta_D(a_i,a_{i+1})}
+\sum_{s+1\le i\le t-2\atop i\notin U_1\cup U_2}
z^{1-\delta_D(a_i,a_{i+1})}.
\relabel{le3-1-eq1}
\end{eqnarray}
(ii) If $m\in \Z^+\setminus V$ and $m>y$ for all $y\in V\setminus R_D[a]$, then 
\begin{eqnarray}
\Delta_{D_{a\rightarrow m},\pi_0}(z)=
z^{2-\delta_D(a_{t-1},a_{t})}+
\sum_{s\le i\le t-2}z^{1-\delta_D(a_i,a_{i+1})}.
\relabel{le3-1-eq2}
\end{eqnarray}
\end{lem}

\proof (i). We will prove this result 
by applying (\ref{def-delta-2}).
Note that $\delta_D(a_s,a)=\delta_D(a,a_t)=1$
as $a_s\rightarrow a$ and $a\rightarrow a_t$ in $D$.
For any $i$ with $s+1\le i\le t-2$, by  (\ref{le3-1-beq1}), 
we have 
\begin{equation}
\delta_D(a_i,a)+\delta_D(a,a_{i+1})=
\left \{
\begin{array}{ll}
0,\qquad &\mbox{if }i\in U_1;\\
2,       &\mbox{if }i\in U_2;\\
1,       &\mbox{otherwise}.\\
\end{array}
\right.\relabel{le3-1-eq3}
\end{equation}
Thus (\ref{le3-1-eq1}) follows from 
(\ref{def-delta-2}).

(ii). 
Recall that $F_D(a)=\{v: (a,v)\in A\}$.
By the assumption on $t$,    
$F_D(a)\subseteq \{a_j: t\le j\le n-1\}$.
As $\pi_0=(a_1,a_2,\cdots,a_{n-1})\in \OP(D-a)$,
we have $R_D(a)\subseteq \{a_j: t\le j\le n-1\}$.
Thus $V(D)\setminus R_D[a]\subseteq \{a_j: 1\le j\le t-1\}$. 
By the assumption on $m$,    
$m>a_i$ holds for all $i:1\le i\le t-1$,
implying that 
\begin{equation}
\delta_{D_{a\rightarrow m}}(a_i,m)
+\delta_{D_{a\rightarrow m}}(m,a_{i+1})
=
\left \{
\begin{array}{ll}
1,\qquad &\mbox{if }s\le i\le t-2;\\
2,       &\mbox{if }i=t-1.\\
\end{array}
\right.
\relabel{le3-1-eq4}
\end{equation}
As $\delta_{D_{a\rightarrow m}}(a_i,a_{i+1})
=\delta_{D}(a_i,a_{i+1})$, 
(\ref{le3-1-eq2}) follows from 
(\ref{def-delta-2}) by replacing $D$ by 
$D_{a\rightarrow m}$.
\proofend

Let 
\begin{equation} 
\left \{ 
\begin{array}{l}
Q(a,\pi_0)
=\{s+1\le i\le t-2: a_i>a>a_{i+1}, 
(a_i,a_{i+1})\in A\};\\
p(a,\pi_0)=(1-\delta_D(a_{s},a_{s+1}))\delta_D(a,a_{s+1})
-(1-\delta_D(a_{t-1},a_{t}))\delta_D(a,a_{t-1}).
\end{array}
\right. \relabel{def-p}
\end{equation}
When there is no confusion, 
$Q(a,\pi_0)$ and $p(a,\pi_0)$ are simply written as $Q$ 
and $p$ respectively. 
Applying Lemma~\ref{le3-1}, 
we can express 
$\Delta_{D,\pi_0}(z)-\Delta_{D_{a\rightarrow m},\pi_0}(z)$
in terms of $Q$ and $p$.

\begin{pro}\relabel{pro3-1}
If $m\in \Z^+\setminus V$ and 
$m>y$ holds for all $y\in V\setminus R_D[a]$,
then 
$$
\Delta_{D,\pi_0}(z)-\Delta_{D_{a\rightarrow m},\pi_0}(z)
=\left (p+|Q|z^{-1}\right )(z-1)^2.
$$
\end{pro}

\proof By (\ref{le3-1-eq1}) and  (\ref{le3-1-eq2})
in Lemma~\ref{le3-1},
\begin{eqnarray}
& &\Delta_{D,\pi_0}(z)-\Delta_{D_{a\rightarrow m},\pi_0}(z) \nonumber \\
&=&z^{1-\delta_D(a_s,a_{s+1})}(z^{\delta_D(a,a_{s+1})}-1)
+z^{1-\delta_D(a_{t-1},a_{t})}(z^{\delta_D(a_{t-1},a)}-z)
\nonumber \\ & &
+(z^{-1}-1)\sum_{i\in U_1}z^{1-\delta_D(a_i,a_{i+1})}
+(z-1)\sum_{i\in U_2}z^{1-\delta_D(a_i,a_{i+1})}
\nonumber \\
&=&z^{1-\delta_D(a_s,a_{s+1})}(z^{\delta_D(a,a_{s+1})}-1)
+z^{1-\delta_D(a_{t-1},a_{t})}(z^{\delta_D(a_{t-1},a)}-z)
\nonumber \\
& &+|Q|(z^{-1}-1)+(|U_1|-|Q|)(1-z)
+|U_2|(z-1)\nonumber \\
&=&z^{1-\delta_D(a_s,a_{s+1})}(z^{\delta_D(a,a_{s+1})}-1)
+z^{1-\delta_D(a_{t-1},a_{t})}(z^{\delta_D(a_{t-1},a)}-z)
\nonumber \\
& &
+(|U_2|-|U_1|)(z-1)+|Q|z^{-1}(z-1)^2,
\relabel{pro3-1-eq1}
\end{eqnarray}
where the second 
last equality follows from the fact that 
for any $i$ with $s+1\le i\le t-2$, 
$$
\delta_D(a_i,a_{i+1})=
\left \{
\begin{array}{ll}
1,\qquad &\mbox{if }i\in Q\cup U_2;\\
0,&\mbox{if }i\in U_1\setminus Q.
\end{array}
\right.
$$
By definitions of $U_1$ and $U_2$, it can be verified that 
\begin{equation}
|U_2|-|U_1|
=
\left \{ 
\begin{array}{ll}
0,\qquad &\mbox{if }a>a_{s+1} \mbox{ and } a>a_{t-1}; \\
1,       &\mbox{if }a>a_{s+1}\mbox{ and } a<a_{t-1}; \\
-1,       &\mbox{if }a<a_{s+1}\mbox{ and } a>a_{t-1}; \\
0,       &\mbox{if }a<a_{s+1}\mbox{ and } a<a_{t-1}.
\end{array}
\right.
\relabel{pro3-1-eq3}
\end{equation}
Then, by (\ref{pro3-1-eq3}),
\begin{eqnarray}
& &z^{1-\delta_D(a_s,a_{s+1})}(z^{\delta_D(a,a_{s+1})}-1)
+z^{1-\delta_D(a_{t-1},a_{t})}(z^{\delta_D(a_{t-1},a)}-z)
+(|U_2|-|U_1|)(z-1)
\nonumber \\
&=&
\left \{ 
\begin{array}{ll}
0, &\mbox{if }a>a_{s+1} \mbox{ and } a>a_{t-1}; \\
z^{1-\delta_D(a_{t-1},a_{t})}(1-z)+(z-1),       
&\mbox{if }a>a_{s+1}\mbox{ and } a<a_{t-1}; \\
z^{1-\delta_D(a_{s},a_{s+1})}(z-1)-(z-1),     
&\mbox{if }a<a_{s+1}\mbox{ and } a>a_{t-1}; \\
z^{1-\delta_D(a_{s},a_{s+1})}(z-1)
+z^{1-\delta_D(a_{t-1},a_{t})}(1-z),     
\quad &\mbox{if }a<a_{s+1}\mbox{ and } a<a_{t-1}
\end{array}
\right. \nonumber \\
&=&
\left \{ 
\begin{array}{ll}
0, &\mbox{if }a>a_{s+1} \mbox{ and } a>a_{t-1}; \\ 
(\delta_D(a_{t-1},a_{t})-1)(z-1)^2,      
&\mbox{if }a>a_{s+1}\mbox{ and } a<a_{t-1}; \\
(1-\delta_D(a_{s},a_{s+1}))(z-1)^2,    
&\mbox{if }a<a_{s+1}\mbox{ and } a>a_{t-1}; \\       
(\delta_D(a_{t-1},a_{t})-\delta_D(a_{s},a_{s+1}))(z-1)^2,
\quad &\mbox{if }a<a_{s+1}\mbox{ and } a<a_{t-1}.
\end{array}
\right.
\relabel{pro3-1-eq4}
\end{eqnarray}
By 
(\ref{def-p}), (\ref{pro3-1-eq1}) and (\ref{pro3-1-eq4}),
the result holds. 
\proofend

By applying Proposition~\ref{pro3-1},
an expression for 
$\Delta(D,z)-\Delta(D_{a\rightarrow m},z)$
can be obtained. 

\begin{theo}\relabel{th3-1}
If $m\in \Z^+\setminus V$ and $m>y$ holds for all $y\in V\setminus R_D[a]$,
then 
\begin{equation}
\Delta(D,z)-\Delta(D_{a\rightarrow m},z)
=(z-1)^2\sum_{\pi_0\in \OP(D-a)}
\left [p(a,\pi_0)+|Q(a,\pi_0)|z^{-1}\right ]
z^{\delta_{D-a}(\pi_0)}.
\end{equation}\relabel{th3-1-eq}
\end{theo}

\proof Observe that 
\begin{eqnarray}
& &\Delta(D,z)-\Delta(D_{a\rightarrow m},z)
\nonumber \\
&=&\sum_{\pi_1\in \OP(D)}z^{\delta_D(\pi_1)}
-\sum_{\pi_2\in \OP(D_{a\rightarrow m})}
z^{\delta_{D_{a\rightarrow m}}(\pi_2)}
\nonumber \\
&=&\sum_{\pi_0\in \OP(D-a)}
\sum_{\pi_1\in \OP(D,\pi_0)}z^{\delta_D(\pi_1)}
-\sum_{\pi_0\in \OP(D-a)}
\sum_{\pi_2\in \OP(D_{a\rightarrow m},\pi_0)}
z^{\delta_{D_{a\rightarrow m}}(\pi_2)}
\nonumber \\
&=&\sum_{\pi_0\in \OP(D-a)}
z^{\delta_{D-a}(\pi_0)}
\Delta_{D,\pi_0}(z)
-\sum_{\pi_0\in \OP(D-a)}
z^{\delta_{D-a}(\pi_0)}
\Delta_{D_{a\rightarrow m},\pi_0}(z)
\nonumber \\
&=& 
\sum_{\pi_0\in \OP(D-a)}
\left [\Delta_{D,\pi_0}(z)
-\Delta_{D_{a\rightarrow m},\pi_0}(z)
\right ]z^{\delta_{D-a}(\pi_0)} 
\nonumber \\
&=&
(z-1)^2\sum_{\pi_0\in \OP(D-a)}
\left [p(a,\pi_0)+|Q(a,\pi_0)|z^{-1}\right ]
z^{\delta_{D-a}(\pi_0)},
\relabel{cor3-1-eq1}
\end{eqnarray}
where the last equality follows from 
Proposition~\ref{pro3-1}.
\proofend

\subsection{Compare 
$D$ with $D_{a\rightarrow \repa}$ for some $\repa>a$\relabel{sec3-2}
}

Let $D=(V,A)$ be an acyclic digraph with $|V|=n$. 
Recall that for $u\in V(D)$, 
$F_D(u)=\{v\in V: (u,v)\in A\}$.
Let $B_D(u)=\{v\in V: (v,u)\in A\}$ and
$B_D[u]=B_D(u)\cup \{u\}$.
Thus $u$ is a sink of $D$ \iff $F_D(u)=\emptyset$,
and $u$ is a source of $D$ \iff $B_D(u)=\emptyset$.

A vertex $u$ of $D$ is called 
a {\it turning vertex} 
if either $F_D(u)=\emptyset$ 
or $\min F_D(u)\ge 2+\max (\P_D(u))$ holds,
where 
\begin{equation}\relabel{sec4-eq1}
\P_D(u)=B_D[u]  \cup 
\{c\in V: \exists b<c, (c,b)\in A\}.
\end{equation}

In this subsection, we always assume that $a$ is a 
turning vertex of $D$ and 
$\repa$ is a number in $\Z^+\setminus V$ such that
$\repa>\max \P_D(a)$ whenever $F_D(a)=\emptyset$,
and $\min F_D(a)>\repa>\max \P_D(a)$
otherwise. 
Thus $y_1>\repa>y_2$ holds for all $y_1\in F_D(a)$ 
and $y_2\in \P_D(a)$. 
Clearly $\repa>a$ holds, as $a\in B_D[a]\subseteq \P_D(a)$.
In this section, 
the assumptions on $a$ and $\repa$ will not be mentioned 
again and 
we shall compare $D$ with $D_{a\rightarrow \repa}$
under this assumption. 

For $u\in V$,  let $\W(D,u)=\{W\in \W(D): u\in W\}$.
So $\W(D,u)=\W(D)\setminus \W(D-u)$, and 
$\W(D,u)=\emptyset$ iff $\W(D)=\W(D-u)$. 

\begin{lem}\relabel{lem4-1}
$\W(D_{a\rightarrow \repa},\repa)=\emptyset$ and so 
$\W(D_{a\rightarrow \repa})=\W(D-a)$.
\end{lem}

\proof Clearly 
$\W(D_{a\rightarrow \repa}) 
=\W(D-a) \cup \W(D_{a\rightarrow \repa},\repa)$.  
Thus it suffices to prove that 
$\W(D_{a\rightarrow \repa},\repa)=\emptyset$, 
i.e., 
$\repa\not\in W$ for every $W\in \W(D_{a\rightarrow \repa})$.

Suppose that $W=\{\repa,b,c\}\in \W(D_{a\rightarrow \repa})$,
where $b<c$. 
Assume that $\repa=\max W$.
Then 
$\repa\rightarrow b$ in $D_{a\rightarrow \repa}$
by definition of $\W(D_{a\rightarrow \repa})$.   
But $\repa\rightarrow b$ in $D_{a\rightarrow \repa}$ implies that 
$a\rightarrow b$ in $D$ and so $b\in F_D(a)$.
By the given condition on $\repa$, 
$\repa<\min F_D(a) \le b$, 
contradicting the assumption that $\repa=\max W>b$.
Hence $\repa<\max W$ and so $\max W=c$. 

If $\repa=\min W$, then, by definition of 
$\W(D_{a\rightarrow \repa})$, 
$c>b>\repa$ and $c\rightarrow \repa$ in $D_{a\rightarrow \repa}$,
where the later implies that 
 $c\rightarrow a$ in $D$. 
 So $c\in B_D(a)\subseteq \P_D(a)$.
 By the given condition on $\repa$, we have 
 $\repa>\max \P_D(a)\ge c$, a contradiction.
 
 By the above conclusions, we have 
$\min W<\repa<\max W$, i.e.,  $b<\repa<c$.
 As $W\in \W(D_{a\rightarrow \repa})$, we have 
 $c\rightarrow b$ in both $D_{a\rightarrow \repa}$ and 
 $D$.
Thus  $c\in \P_D(a)$.
 But $\repa>\max \P_D(a)$ implies that $\repa>c$,
 a contradiction again. 
 
 Hence the result holds.
\proofend

\begin{lem}\relabel{lem4-2} 
Let $\pi_0=(a_1,a_2,\cdots,a_{n-1})\in \OP(D-a)$
and $s$ and $t$ be the numbers defined in 
Subsubsection~\ref{sec3-1-1} with respect to $a$ and $\pi_0\in \OP(D-a)$.
Then 
\begin{enumerate}
\item\relabel{lem4-2-no1} 
$Q(\repa,\pi_0)=\emptyset$;

\item\relabel{lem4-2-no2} 
$p(a,\pi_0)-p(\repa,\pi_0)=1$ if $\{a,a_{s+1},a_s\}\in \W(D)$,
and $p(a,\pi_0)-p(\repa,\pi_0)=0$ otherwise. 
\end{enumerate}
\end{lem}

\proof 
\ref{lem4-2-no1}
By definition, 
$$
Q(\repa,\pi_0)
=\{s+1\le i\le t-2: a_i>\repa>a_{i+1},a_i\rightarrow a_{i+1}\}.
$$
Assume that $k\in Q(\repa,\pi_0)$.
Then $a_k>a_{k+1}$ and $a_k\rightarrow a_{k+1}$,
implying that $a_k\in \P_D(a)$.
By the assumption on $\repa$, we have 
$\repa>\max \P_D(a)\ge a_k$.
However, $k\in Q(\repa,\pi_0)$ implies that 
$\repa<a_k$, a contradiction. 
Thus $Q(\repa,\pi_0)=\emptyset$.

\ref{lem4-2-no2}
By definition of $p(a,\pi_0)$, we have 
\begin{eqnarray}
p(a,\pi_0)-p(\repa,\pi_0)
&=&(1-\delta_D(a_s,a_{s+1}))
\left [\delta_D(a,a_{s+1})
-\delta_{D_{a\rightarrow \repa}}(\repa,a_{s+1})\right ]
\nonumber \\
& &+(1-\delta_D(a_{t-1},a_{t}))
\left [\delta_{D_{a\rightarrow \repa}}(\repa,a_{t-1})
-\delta_D(a,a_{t-1})\right ].
\relabel{lem4-2-eq3}
\end{eqnarray}

\noindent {\bf Claim 1}: 
$p(a,\pi_0)-p(\repa,\pi_0)
=(1-\delta_D(a_s,a_{s+1}))
\left [\delta_D(a,a_{s+1})
-\delta_{D_{a\rightarrow \repa}}(\repa,a_{s+1})\right ]$.

By (\ref{lem4-2-eq3}), it suffices to show that 
$(1-\delta_D(a_{t-1},a_{t}))
\left [\delta_{D_{a\rightarrow \repa}}(\repa,a_{t-1})-\delta_D(a,a_{t-1})\right ]=0.
$
Suppose that it does not hold.
Then  $\delta_D(a_{t-1},a_{t})=0$.
Thus  $t<n$ and 
$a_{t-1}>a_t$. 
By the assumption on $t$, we have $(a,a_t)\in A$,
implying that $a_t\in F_D(a)$.
Since $a<\repa<\min F_D(a)$,   
we have $a<\repa<a_t$.
As $a_t<a_{t-1}$, we have $a<\repa<a_{t-1}$ and 
$$
\delta_{D_{a\rightarrow \repa}}(\repa,a_{t-1})=\delta_D(a,a_{t-1})=1.
$$
So $(1-\delta_D(a_{t-1},a_{t}))
\left [\delta_{D_{a\rightarrow \repa}}(\repa,a_{t-1})
-\delta_D(a,a_{t-1})\right ]=0$,
a contradiction. 
Hence Claim 1 holds.

\noindent {\bf Claim 2}:
$p(a,\pi_0)-p(\repa,\pi_0)\ge 0$. 

As $\repa>a$,
$\delta_D(a,a_{s+1})-\delta_{D_{a\rightarrow \repa}}(\repa,a_{s+1})\ge 0$.
Then Claim 2 follows from Claim 1.

\noindent {\bf Claim 3}: 
 $p(a,\pi_0)-p(\repa,\pi_0)=1$
\iff $\{a,a_s,a_{s+1}\}\in \W(D)$.

By Claims 1 and 2, $p(a,\pi_0)-p(\repa,\pi_0)\in \{0,1\}$.

Assume that
$\{a,a_s,a_{s+1}\}\in \W(D)$.
By definition of $s$, $a_s\rightarrow a$ in $D$.
By definition of $\W(D)$, $a_s>a_{s+1}>a$, 
$a\not\rightarrow a_{s+1}$ and 
$a_s\not\rightarrow a_{s+1}$
in $D$.
So $\delta_D(a_s,a_{s+1})=0$ and $\delta_D(a,a_{s+1})=1$.
As $a_s\in B_D[a]\subseteq \P_D(a)$, 
by the assumption on $\repa$, $\repa>a_s$ holds,
implying that $\repa>a_s>a_{s+1}$.
Since $a\not\rightarrow a_{s+1}$ in $D$,
we have $\repa\not\rightarrow a_{s+1}$ in $D_{a\rightarrow \repa}$.
Thus
$\delta_{D_{a\rightarrow \repa}}(\repa,a_{s+1})=0$.
By Claim 1, we have $p(a,\pi_0)-p(\repa,\pi_0)=1$ 

Now assume that $p(a,\pi_0)-p(\repa,\pi_0)=1$.
By Claim 1,
$\delta_D(a_s,a_{s+1})=0$ and 
$\delta_D(a,a_{s+1})
-\delta_{D_{a\rightarrow \repa}}(\repa,a_{s+1})=1$,
where the later implies that  $\delta_D(a,a_{s+1})=1$.
Observe that $\delta_D(a_s,a_{s+1})=0$ implies  that 
$a_s>a_{s+1}$ and $a_s\not \rightarrow a_{s+1}$,
and $\delta_D(a,a_{s+1})=1$ implies that  
$a<a_{s+1}$ or $a\rightarrow a_{s+1}$ in $D$.
However, if $a\rightarrow a_{s+1}$ in $D$, then 
$\repa\rightarrow a_{s+1}$ in $D_{a\rightarrow \repa}$, 
implying that 
$\delta_D(a,a_{s+1})
-\delta_{D_{a\rightarrow \repa}}(\repa,a_{s+1})=1-1=0$,
a contradiction. Thus 
$a_s>a_{s+1}>a$, but $a_s\not \rightarrow a_{s+1}$ 
and $a\not \rightarrow a_{s+1}$
in $D$.
By definition of $s$,  $a_s\rightarrow a$.
Hence  
$\{a,a_s,a_{s+1}\}\in \W(D)$ and the claim holds.
\proofend

For an integer $j$ with $0\le j\le n-1$, 
\begin{enumerate}
\item let $c_j(D,a)$ be the number of 
$\pi=(a_1,\cdots,a_i,a,a_{i+1},\cdots,a_{n-1})\in \OP(D)$ 
such that $\delta_{D}(\pi)=j$ and 
$\{a,a_i,a_{i+1}\}\in \W(D)$ for some $i$ with $1\le i\le n-1$, 
where 
$(a_i,a)\in A$;

\item let $c'_j(D,a)$ be the number of 
$\pi=(a_1,\cdots,a_i,a,a_{i+1},\cdots,a_{n-1})$ 
such that $\delta_{D}(\pi)=j$ and 
$\{a,a_i,a_{i+1}\}\in \W(D)$ for some $i$ with $1\le i\le n-1$, 
where 
$(a_i,a_{i+1})\in A$.
\end{enumerate}

Clearly $c_j(D,a)+C'_j(D,a)$ is not more than 
the number of $\pi$'s in $\OP(D)$ with $\delta_D(\pi)=j$, and 
$c_j(D,a)=C'_j(D,a)=0$ whenever $\W(D,a)=0$.

\begin{lem}\relabel{le3-4}
$c_{j}(D,a)=0$ for $j=0,1$,  
and $c'_j(D,a)=0$ for $j\ge n-2$.
\end{lem}

\proof For any 
$\pi=(a_1,\cdots,a_i,a,a_{i+1},\cdots,a_{n-1})\in \OP(D)$, 
if $\{a,a_i,a_{i+1}\}$ is a member 
of $\W(D)$ with $a<a_{i+1}<a_i$  
and $(a_{i},a)\in A$, 
then $\delta_D(a_i,a)=\delta_D(a,a_{i+1})=1$,
implying that $\delta_D(\pi)\ge 2$. 
Thus $c_{j}(D,a)=0$ for $j\le 1$
by definition of $c_j(D,a)$.

For any $\pi=(a_1,\cdots,a_i,a,a_{i+1},\cdots,a_{n-1})
\in \OP(D)$, 
if $\{a,a_i,a_{i+1}\}$ is a member 
of $\W(D)$ with $a_{i+1}<a<a_i$  
and $(a_{i},a_{i+1})\in A$, 
then $a_i\not \approx a$ and $a\not\approx a_{i+1}$,
implying that 
$\delta_D(a_i,a)=\delta_D(a,a_{i+1})=0$.
Thus $\delta_D(\pi)\le n-3$.
By definition of $c'_j(D,a)$,
$c'_{j}(D,a)=0$ for $j\ge n-2$.
\proofend

\begin{theo}\relabel{th4-1}
Assume that $n=|V|\ge 3$. Then
\begin{equation}
\Delta(D,z)-\Delta(D_{a\rightarrow \repa},z)
=(z-1)^2\sum_{0\le j\le n-3} (c_{j+2}(D,a)+c'_{j}(D,a)) z^j.
\relabel{th4-1-eq1}
\end{equation}
Furthermore, $\Delta(D,z)=\Delta(D_{a\rightarrow \repa},z)$
\iff $\W(D,a)=\emptyset$.
\end{theo}

\proof 
Let $m$ be a number in $\Z^+\setminus V$ such that $m>y$ for all 
$y\in V\setminus R_D[a]$.
By Theorem~\ref{th3-1}, we have 
\begin{equation}
\Delta(D,z)-\Delta(D_{a\rightarrow m},z)
=(z-1)^2\sum_{\pi_0\in \OP(D-a)}
\left [
p(a,\pi_0)+|Q(a,\pi_0)|z^{-1}
\right] z^{\delta_{D-a}(\pi_0)}.
\relabel{th4-1-eq0-1}
\end{equation}
By Lemma~\ref{lem4-2}~\ref{lem4-2-no1}, 
$Q(\repa,\pi_0)=\emptyset$.
Replacing 
$D$ by $D_{a\rightarrow \repa}$ in (\ref{th4-1-eq0-1}) 
gives that 
\begin{equation}
\Delta(D_{a\rightarrow \repa},z)-\Delta(D_{a\rightarrow m},z)
=(z-1)^2\sum_{\pi_0\in \OP(D-a)}
p(\repa,\pi_0) z^{\delta_{D-a}(\pi_0)}.
\relabel{th4-1-eq0-2}
\end{equation}
By (\ref{th4-1-eq0-1}) and (\ref{th4-1-eq0-1}),
$\Delta(D,z)-\Delta(D_{a\rightarrow \repa},z)$ 
has the following expression:
\begin{equation}
(z-1)^2\sum_{\pi_0\in \OP(D-a)}
\left [
p(a,\pi_0)-p(\repa,\pi_0)+|Q(a,\pi_0)|z^{-1}
\right] z^{\delta_{D-a}(\pi_0)}.
\relabel{th4-1-eq2}
\end{equation}

The proof will be completed 
by establishing the following claims.    

\inclaim For each $\pi_0=(a_1,a_2,\cdots,a_{n-1})\in \OP(D-a)$, 
$p(a,\pi_0)-p(\repa,\pi_0)\in \{0,1\}$,
and $p(a,\pi_0)-p(\repa,\pi_0)=1$ \iff 
$\{a,a_s,a_{s+1}\}\in \W(D)$, where $(a_s,a)\in A$.

Claim~\thecountclaim\ follows from 
Lemma~\ref{lem4-2}~\ref{lem4-2-no2}.

\inclaim 
$\sum\limits_{\pi_0\in\OP(D-a)}
(p(a,\pi_0)-p(\repa,\pi_0)) 
z^{\delta_{D-a}(\pi_0)}
=\sum\limits_{j=0}^{n-3} c_{j+2}(D,a)x^j$.

Let $\OP^*(D-a)$ be the set of those $\pi_0\in \OP(D-a)$ 
with $p(a,\pi_0)-p(\repa,\pi_0)=1$,
and let $q_j$ be the number of $\pi_0$'s in $\OP^*(D-a)$
with $\delta_{D-a}(\pi_0)=j$, where $0\le j\le n-2$.
Then, by Claim 1, 
\begin{align}
\sum_{\pi_0\in\OP(D-a)}(p(a,\pi_0)-p(\repa,\pi_0)) 
z^{\delta_{D-a}(\pi_0)}
&=\sum_{j=0}^{n-2}
\sum_{\pi_0\in \OP^*(D-a)
\atop \delta_{D-a}(\pi_0)=j}
z^{j}
=\sum_{j=0}^{n-2} q_jz^j.
\relabel{th4-1-eq5}
\end{align}
For each $\pi_0\in \OP^*(D-a)$ with $\delta_{D-a}(\pi_0)=j$, 
$\pi=(a_1,\cdots,a_s,a,a_{s+1},\cdots,a_{n-1})$
is a member of $\OP(D)$.  
By Claim 1, $\{a,a_s,a_{s+1}\}\in \W(D)$ 
with 
$(a_s,a)\in A$.
Thus, $\delta_D(a_s,a)=\delta(a,a_{s+1})=1$ 
but $\delta_{D-a}(a_s,a_{s+1})=0$,
implying that 
$\delta_D(\pi)=\delta_{D-a}(\pi_0)+2=j+2$.
As $\delta_{D-a}(a_s,a_{s+1})=0$, we have  
$\delta_{D-a}(\pi_0)\le n-3$ and so $q_{n-2}=0$.
On the other hand, for any $\pi=(a_1,\cdots,a_i,a,a_{i+1},\cdots,a_{n-1})\in \OP(D)$
with $\delta_D(\pi)=j+2$
and $\{a_i,a,a_{i+1}\}\in \W(D)$, 
where $a_i>a_{i+1}>a$, 
$\pi_0=(a_1,\cdots,a_i,a_{i+1},\cdots,a_{n-1})$
is a $(D-a)$-respecting ordering with $s=i$
and $\delta_{D-a}(\pi)=j$.
Thus, by definition, $q_j=c_{j+2}(D,a)$ holds
and so Claim~\thecountclaim\ holds.

\inclaim 
$\sum\limits_{\pi_0\in \OP(D-a)}
|Q(a,\pi_0)|z^{\delta_{D-a}(\pi_0)-1}
=\sum\limits_{j=0}^{n-3} c'_j(D,a)z^{j}$.

By definition, for each 
$\pi_0=(a_1,a_2,\cdots,a_{n-1})\in \OP(D-a)$,
$|Q(a,\pi_0)|$ is the number of integers $i$  
with $s+1\le i\le t-2$
such that $a_i>a>a_{i+1}$ and $(a_i,a_{i+1})\in A$.
As $(a_i,a_{i+1})\in A$, we have $\delta_{D-a}(\pi_0)\ge 1$.
As $s+1\le i\le t-2$,  
the definitions of $s$ and $t$ imply that 
$D[\{a_i,a,a_{i+1}\}]$ has 
only one arc, i.e., $(a_i,a_{i+1})$.
Thus $\{a_i,a,a_{i+1}\}\in \W(D)$.
Clearly, $Q(a,\pi_0)>0$ implies that  
$\delta_{D-a}(\pi_0)\ge 1$.
Thus,
\begin{align}
&\sum_{\pi_0\in \OP(D-a)}
|Q(a,\pi_0)|z^{\delta_{D-a}(\pi_0)-1}
=\sum_{j=0}^{n-2} 
\sum_{\pi_0\in \OP(D-a)\atop \delta_{D-a}(\pi_0)=j} 
Q(a,\pi_0)z^{j-1}
=\sum_{j=0}^{n-3} q'_jz^j,
\relabel{th4-1-eq6}
\end{align}
where $q'_j$ is the number of order pairs 
$(\pi_0,i)$, where $\pi_0\in \OP(D-a)$ 
with $\delta_{D-a}(\pi_0)=j+1$
and $i$ is an integer with $s+1\le i\le t-2$
such that $\{a_i,a,a_{i+1}\}\in \W(D)$, where  
$(a_i,a_{i+1})\in A$.

For each $\pi_0=(a_1,a_2,\cdots,a_{n-1})\in \OP(D-a)$
with $\delta_{D-a}(\pi_0)=j+1$,
if $i$ is an integer with $s+1\le i\le t-2$
such that $\{a_i,a,a_{i+1}\}\in \W(D)$, where  
$(a_i,a_{i+1})\in A$,
then $\pi=(a_1,\cdots,a_i,a,a_{i+1},\cdots,a_{n-1})$ 
is a member of $\OP(D)$.  
As $\delta_D(a_i,a)=\delta_D(a,a_{i+1})=0$
but $\delta_{D-a}(a_i,a_{i+1})=1$, 
we have $\delta_{D}(\pi)=\delta_{D-a}(\pi_0)-1=j$.

On the other hand, for each 
$\pi=(a_1,\cdots,a_i,a,a_{i+1},\cdots,a_{n-1})$,
if $\{a_i,a,a_{i+1}\}\in \W(D)$, where 
$(a_i,a_{i+1})\in A$, 
by definitions of $s$ and $t$, 
we have $s+1\le i\le t-2$
and $\pi_0=(a_1,\cdots,a_i,a_{i+1},\cdots,a_{n-1})$
is a member of $\OP(D-a)$. 
As $\delta_D(a_i,a)=\delta_D(a,a_{i+1})=0$
and $\delta_{D-a}(a_i,a_{i+1})=1$, 
we have $\delta_{D}(\pi)=\delta_{D-a}(\pi_0)-1=j$
whenever $\delta_{D-a}(\pi_0)=j+1$.


By the assumption on $q'_j$ and the above arguments, $q'_j$ equals 
the number of members  
$\pi=(a_1,\cdots,a_i,a,a_{i+1},\cdots,a_{n-1})$ of $\OP(D)$  
with $\delta_D(\pi)=j$
such that $\{a_i,a,a_{i+1}\}\in \W(D)$, 
where 
$(a_i,a_{i+1})\in A$.
By definition of $c'_j(D,a)$,
we have $q'_j=c'_j(D,a)$.
Then, by (\ref{th4-1-eq6}), Claim~\thecountclaim\ holds.

By (\ref{th4-1-eq2}) and Claims 2 and 3,
(\ref{th4-1-eq1}) holds.

\inclaim If $\W(D,a)\ne \emptyset$,
then $c_{j+2}(D,a)+c'_{j}(D,a)>0$ for some $j$.

Assume that $W=\{a,b,c\}\in\W(D)$, where $b<c$.
If $a>c$, then $(a,b)\in A$, implying that $b\in F_D(a)$.
But  $a$ is a turning vertex of $D$, implying that 
$a<y$ for all $y\in F_D(a)$, a contradiction.
Thus, either $a<b<c$ or $b<a<c$.

Suppose that $a<b<c$. As $\{a,b,c\}\in W(D)$, by definition 
of $\W(D)$, $(c,a)\in A$ and 
$b\not \approx_D  a$  and $b\not \approx_D  c$.
It is easy to check that there exists $\pi=(a_1,\cdots,a_s,a,a_{s+1},\cdots,a_{n-1})\in \OP(D)$,
where $a_s=c$ and $a_{s+1}=b$.
Thus $\{a,a_s,a_{s+1}\}\in \W(D)$.
Let $\pi_0=\pi-a$, i.e., 
$\pi_0=(a_1,\cdots,a_s,a_{s+1},\cdots,a_{n-1})$.
Clearly $\pi_0\in \OP(D-a)$,
$\delta_{D}(\pi)\ge \delta_{D}(a_s,a)
+\delta_{D}(a,a_{s+1})=2$ 
and $\delta_D(\pi)=\delta_{D-a}(\pi_0)+2\ge 2$.
By definition, 
$c_{j+2}(D,a)>0$ for some $j$ with $0\le j\le n-3$.

Now suppose that $b<a<c$. 
As $\{a,b,c\}\in W(D)$, by definition 
of $\W(D)$, $(c,b)\in A$ and 
$a\not \approx_D  b$  and $a\not \approx_D  c$.
It is easy to check that there exists $\pi=(a_1,\cdots,a_i,a,a_{i+1},\cdots,a_{n-1})\in \OP(D)$,
where $a_i=c$ and $a_{i+1}=b$.
Let $\pi_0=(a_1,\cdots,a_i,a_{i+1},\cdots,a_{n-1})$.
Clearly $\pi_0\in \OP(D-a)$ and 
$\delta_D(\pi)=\delta_{D-a}(\pi_0)-1\le n-3$.
By definition, 
$c'_j(D,a)>0$ for some $j$ with $0\le j\le n-3$.

Thus Claim~\thecountclaim\ holds.    
If $\W(D,a)=\emptyset$, by definition of $c_j(D,a)$
and $c'_j(D,a)$, we have 
$c_j(D,a)=c'_j(D,a)=0$
for all $i=0,1,\cdots,n-1$. 
By this fact and Claim 4, 
$\W(D,a)=\emptyset$ 
\iff $\Delta(D,z)=\Delta(D_{a\rightarrow r},z)$.
\proofend

Applying Theorem~\ref{th4-1}
and the following result, we will obtain an expression for 
$\Psi(D,x)-\Psi(D_{a\rightarrow \repa},x)$
in terms of $c_{j+2}(D,a)+c'_{j}(D,a)$ for $j=0,1,\cdots,n-3$.

\begin{lem}\relabel{le4-3}
Let  $D_1$ and $D_2$ be any two acyclic digraphs
of order $n$.
\begin{enumerate}
\item 
If $\Delta(D_1,z)-\Delta(D_2,z)=
t_0+t_1z+\cdots+t_{n-1}z^{n-1}$,
then 
\begin{equation}\relabel{le4-3-eq1}
\Psi(D_1,x)-\Psi(D_2,x)=\sum_{i=0}^{n-1} t_i 
{x+i\choose n};
\end{equation}
\item 
if $\Delta(D_1,z)-\Delta(D_2,z)=(z-1)^2P(z)$,
where $P(z)=d_0+d_1z+\cdots+d_{n-3}z^{n-3}$,
then 
\begin{equation}\relabel{le4-3-eq2}
\Psi(D_1,x)-\Psi(D_2,x)=\sum_{i=0}^{n-3} d_i 
{x+i\choose n-2}.
\end{equation}
\end{enumerate}
\end{lem}

\proof (i). Assume that 
$$
\Delta(D_2,z)=\sum_{i=0}^{n-1} b_iz^{i}.
$$
Then, by the given condition,
$$
\Delta(D_1,z)=\sum_{i=0}^{n-1} (b_i+t_i)z^{i}.
$$
By the relation between $\Delta(D_i,z)$ and $\Psi(D_i,x)$,
we have 
$$
\Psi(D_1,x)=\sum_{i=0}^{n-1} (b_i+t_i){x+i\choose n},
\quad
\Psi(D_2,x)=\sum_{i=0}^{n-1} b_i{x+i\choose n}.
$$
Thus the result holds.

(ii). Note that 
$$
\Delta(D_1,z)-\Delta(D_2,z)
=(z-1)^2\sum_{i=0}^{n-3}d_iz^i
=\sum_{i=0}^{n-3} 
\left ( d_iz^{i+2}-2d_iz^{i+1}+d_iz^{i}\right ).
$$
Then, the result in (i) implies that 
\begin{eqnarray}
\Psi(D_1,x)-\Psi(D_2,x)
&=&\sum_{i=0}^{n-3} d_i \left [
{x+i+2\choose n}-2{x+i+1\choose n}
+{x+i\choose n}
\right ]
\nonumber \\
&=&\sum_{i=0}^{n-3} d_i {x+i\choose n-2}.
\end{eqnarray}
\proofend

\begin{theo}\relabel{th4-2}
Assume that $n=|V|\ge 3$. Then 
\begin{equation}\relabel{th4-2-eq1}
\Psi(D,x)-\Psi(D_{a\rightarrow \repa},x)
=\sum_{j=0}^{n-3} (c_{j+2}(D,a)+c'_{j}(D,a)){x+j\choose n-2}.
\end{equation}
Furthermore,  
$\Psi(D,x)=\Psi(D_{a\rightarrow \repa},x)$
\iff $\W(D,a)=\emptyset$.
\end{theo}

\proof
The result follows directly from 
Theorem~\ref{th4-1} and Lemma~\ref{le4-3}(ii).
\proofend

Let $D_1=(V_1,A_1)$ be an acyclic digraph and 
$V'\subseteq V_1$.
Let $D_2=(V_2,A_2)$ be an acyclic digraph obtained from 
$D_1$ by relabeling each $u\in V'$ by $\mu(u)$,
where $\mu$ is a bijection from $V'$ to 
$V''$, where $V''$ is some subset of 
$\Z^+\setminus V_1$ with $|V''|=|V'|$.    
Write $D_1\succeq D_2$ if conditions (a) and (b) below 
are satisfied:
\begin{enumerate}
\item[(a)] for any $3$-element subset $W$ of $V_1$, 
if $W\not\in \W(D_1)$, then $W'\not\in \W(D_2)$, 
where $W'=(W\setminus V')\cup \{\mu(u): u\in W\cap V'\}$;

\item[(b)] 
$\Delta(D_1,z)-\Delta(D_2,z)=(z-1)^2P(z)$, 
where $P(z)=0$ or $P(z)$ is a polynomial of degree 
at most $n_1-3$
without negative coefficients, where $n_1=|V_1|$;
furthermore, $P(z)=0$ \iff 
$|\W(D_1)|=|\W(D_2)|$.
\end{enumerate}

\begin{pro}\relabel{prop4-0}
If $D_1\succeq D_2$ and 
$D_2\succeq D_3$, then  $D_1\succeq D_3$.
\end{pro}

\begin{pro}\relabel{prop4-1}
Assume that $D_1\succeq D_2$. Then 
\begin{enumerate}
\item $|\W(D_1)|\ge |\W(D_2)|$;
\item if $|\W(D_1)|=|\W(D_2)|$, 
then $\Delta(D_1,z)=\Delta(D_2,z)$ and 
$\Psi(D_1,x)=\Psi(D_2,x)$;
\item if
$|\W(D_1)|>|\W(D_2)|$, then 
there exists non-negative integers $d_0,d_1,\cdots,d_{n_1-3}$
such  that
$$
\Delta(D_1,z)-\Delta(D_2,z)=(z-1)^2
\sum_{i=0}^{n_1-3} d_iz^i
$$
and 
$$
\Psi(D_1,x)-\Psi(D_2,x)
=\sum_{i=0}^{n_1-3} d_i{x+i\choose n_1-2},
$$ 
where $d_i>0$ for some $i$.
\end{enumerate}
\end{pro}

By applying Lemma~\ref{lem4-1} and Theorem~\ref{th4-1},
we get the following conclusion 
on $D$ and $D_{a\rightarrow \repa}$.

\begin{cor}\relabel{cor4-1}
$D\succeq  D_{a\rightarrow \repa}$.
\end{cor}

\subsection{
Complete the proof of Theorem~\ref{th-new2}
\relabel{sec3-3}
}

Let $D=(V,A)$ be an acyclic digraph.  For $S\subseteq V$,  
$S$ is said to be {\it ideal} in  
$D$ if either $S=\emptyset$ or 
the following conditions  are satisfied:
\begin{enumerate}
\item[(i.1)] for each $y\in S$, $R_D(y)\subseteq S$;
\item[(i.2)] for each $y\in S$, either $F_D(y)=\emptyset$
or $y<\min F_D(y)$; and 
\item[(i.3)] either $S=V$ and $\min V\ge 2$ 
or $\min S\ge 2+\max (V\setminus S)$.
\end{enumerate}


\begin{pro}\relabel{pro5-0}
Let $S\subseteq V$ be ideal in $D$. Then 
$\Re(D)=\Re(D-S)$
and $\W(D)=\W(D-S)$.
\end{pro}

\proof 
We just need to consider the case that 
$S\ne \emptyset$.    
As $S$ is ideal in $D$, 
it is easy to verify that $\Re(D)=\Re(D-S)$.

It is clear that $\W(D-S)\subseteq \W(D)$.
Assume that $W\in \W(D)$ and $W\cap S\ne \emptyset$.
As $\min S\ge 2+\max (V\setminus S)$, 
we have $\max W\in S$.
Let $c=\max W$ and $a=\min W$. So $c>a$.
By definition of $\W(D)$, $(c,a)\in A$
and so $a\in F_D(c)$.
As $S$ is ideal, $c<\min F_D(c)\le a$, a contradiction.
\proofend

Assume that $\min S=+\infty$ whenever $S=\emptyset$.

\begin{pro}\relabel{pro5-1}
Let $S\subseteq V$ be ideal in $D$
and let $u\in V\setminus S$ 
with $F_D(u)\subseteq S$. Then 
\begin{enumerate}
\item\relabel{pro5-1-no2} 
$\P_D(u)\subseteq V\setminus S$ and
$u$ is  a turning vertex of $D$;

\item\relabel{pro5-1-no3} 
if $V=S\cup \{u\}$, then 
$S\cup \{u'\}$ is ideal in $D_{u\rightarrow u'}$
for any $u'\in \Z^+$ with $0<u'<\min S$;
\item\relabel{pro5-1-no4} 
if $V\ne S\cup \{u\}$ and
$\min S\ge 3+\max (V\setminus S)$, then 
$S\cup \{u'\}$ is ideal in $D_{u\rightarrow u'}$
for any $u'\in \Z^+$  with 
$2+\max (V\setminus S) \le u'<\min S$.
\end{enumerate}
\end{pro}

\proof 
\ref{pro5-1-no2}
The result is trivial if $S=\emptyset$.
So we assume that $S\ne \emptyset$.
As $S$ is ideal in $D$ and $u\notin S$, we 
have $B_D(u)\subseteq V\setminus S$.
For any $(c,b)\in A$, if $c\in S$, then 
$b\in S$ by condition (i.1) and so $c<b$ by condition (i.2).
Thus, $(c,b)\in A$ and $c>b$ imply that $c\notin S$.
Therefore, 
\begin{eqnarray}\relabel{pro5-1-eq1}
\P_D(u)&=&B_D[u]\cup 
\{c\in V: \exists b<c, (c,b)\in A\}
\subseteq V\setminus S.
\end{eqnarray}
As $F_D(u)\subseteq S$ and $S$ is ideal in $D$,
we have 
\begin{eqnarray}\relabel{pro5-1-eq2}
\min F_D(u)\ge \min S\ge 
2+\max (V\setminus S)\ge 2+\max P_D(u).
\end{eqnarray}
Thus, $u$ is a turning vertex of $D$.

\ref{pro5-1-no3} 
This is trivial to verify.    

\ref{pro5-1-no4} The result is trivial when $S=\emptyset$.
Now assume that $S\ne \emptyset$.
Let $S'=S\cup \{u'\}$.
By the given condition, to verify if 
$S'$ is ideal in $D_{u\rightarrow u'}$,
it suffices to show that condition (i.3) is satisfied.
As $u'=\min S-1$ and $\min S\ge 3+\max (V\setminus S)$,
we have 
\begin{eqnarray}\relabel{pro5-1-eq3}
\min S'
&=&u'\ge 2+\max (V\setminus S)\ge 
2+\max (V\setminus (S\cup \{u\}))
\nonumber \\
& = &2+\max (V(D_{u\rightarrow u'})\setminus S').
\end{eqnarray}
Thus $S'$ is ideal in $D_{u\rightarrow u'}$.
\proofend

\begin{pro}\relabel{pro5-2}
Let $S\subset V$ be ideal in $D$
and $u$ be a vertex in $V\setminus S$ with $F_D(u)\subseteq S$.
For any $u'\in \Z^+$ with $\max(V\setminus S)<u'<\min S$,
$D\succeq D_{u\rightarrow u'}$ holds.
\end{pro}

\proof By Proposition~\ref{pro5-1}~\ref{pro5-1-no2},
$\P_D(u)\subseteq V\setminus S$ and  
$u$ is a turning vertex of $D$.
Thus, if $\max(V\setminus S)<u'<\min S$,
then $\max \P_D(u)\le \max(V\setminus S)<u'<\min S\le \min F_D(u)$.
Replacing $r$ by $u'$ in Corollary~\ref{cor4-1}
implies that 
$D\succeq D_{u\rightarrow u'}$. 
\proofend

For an acyclic digraph $D=(V,A)$,  
an ordering $\alpha=(u_1,u_2,\cdots,u_n)$ of its vertices
is said to be a {\it sink-elimination ordering}, 
if $u_i$ is a sink of the subdigraph 
$D[V_i]$ for all $i=1,2,\cdots,n-1$,
 where $V_i=\{u_i,u_{i+1},\cdots,u_n\}$.
 Now assume that $\alpha=(u_1,u_2,\cdots,u_n)$  
 is a sink-elimination ordering of $D$
 and $M=n+1+\max V$.
 Let $\Gamma_{D,\alpha}$ denote the sequence 
 $(D_0,D_1,\cdots,D_{n-1})$ of digraphs  
 produced from $D$ according to $\alpha$:
$D_0$ is $D$, and 
for $i=1,2,\cdots,n-1$, 
$D_i$ is the digraph $(D_{i-1})_{u_i\rightarrow M-i}$
(i.e., $D_i$ is 
obtained from $D_{i-1}$ 
by relabeling vertex $u_i$ as $M-i$).
For example, if $D$ is the digraph in Figure~\ref{f2},
then $\alpha=(3,2,4,5,1)$ is a sink-elimination 
ordering of its vertices,
$M=11$ and $\Gamma_{D,\alpha}=(D_0,D_1,\cdots,D_4)$,
where $D_0$ is the digraph in Figure~\ref{f2}, 
$D_1,D_2,D_3$ and $D_4$ 
are shown in Figure~\ref{p3}.

\begin{figure}[htbp]
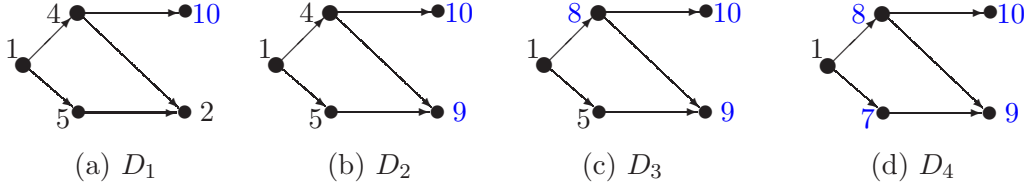

  \centering
\input p3.pic

\vspace{0.3 cm}

(a) $D_1$ \hspace{2 cm} 
(b) $D_2$ \hspace{2 cm} 
(c) $D_3$ \hspace{2.5 cm} 
(d) $D_4$

  \caption{$\Gamma_{D,\alpha}=(D_0,D_1,\cdots,D_4)$
  for $D$ in Figure~\ref{f2} and $\alpha=(3,2,4,5,1)$}
\relabel{p3}
\end{figure}

\begin{theo}\relabel{theo5-2}
Assume that $\Gamma_{D,\alpha}=
(D_0,D_1,\cdots,D_{n-1})$.
Then 
$\Re(D_{n-1})=\emptyset$
and 
$D_i\succeq D_{i+1}$ for all $i=0,1,\cdots,n-2$.
\end{theo}

\proof Let $M=n+1+\max V$. 
By definition,  
$D_i$ is obtained from $D$ 
by relabeling vertex $u_j$ as  $M-j$ for all $j=1,2,\cdots,i$,
where $\alpha=(u_1,u_2,\cdots,u_n)$ is a sink-elimination 
ordering of $D$.
Then $V(D_i)$ is the disjoint union of $S_i$ and $V_{i+1}$,
where $S_i=\{M-j: 1\le j\le i\}$ and 
$V_{i+1}=\{u_j: i+1\le j\le n\}$.

We first prove two claims below.

\noindent {\bf Claim 1}:
$F_{D_i}(u_{i+1})\subseteq S_i$ 
for all $i=0,1,\cdots,n-1$.

As $\alpha$ is a sink-elimination ordering of $D$,
$u_{i+1}$ is a sink of $D[V_{i+1}]$
and so 
$F_D(u_{i+1})\subseteq \{u_1,\cdots,u_{i}\}$.
By definition of $D_i$, we have 
$F_{D_i}(u_{i+1})\subseteq \{M-j: 1\le j\le i \}
=S_i$.
Hence Claim 1 holds.

\noindent {\bf Claim 2}:
$S_i$ is ideal in $D_i$
for all $i=0,1,\cdots,n-1$.

As $S_0=\emptyset$, $S_0$ is ideal in $D_0$.
It is also trivial that $S_1=\{M-1\}$ is ideal in $D_1$,
as $M-1=n+\max V\ge 3+\max V\ge \max (V_2)+3$
and $M-1$ is a sink in $D_1$.

Now assume that $S_{i-1}$ is ideal in $D_{i-1}$, 
where $2\le i\le n-2$. 
We will apply Proposition~\ref{pro5-1}
to show that $S_{i}$ is ideal in $D_{i}$.

Note that $u_i\in V(D_{i-1})\setminus S_{i-1}=V_i$ and 
$D_i$ is obtained from $D_{i-1}$ by relabeling 
$u_i$ as $M-i$. 
Observe that 
$M-i<M-i+1=\min S_{i-1}$ and 
$$
M-i=n+1+\max V -i \ge 3+\max V\ge 3+\max V_i.
$$
By Claim 1,  
$F_{D_{i-1}}(u_i)\subseteq S_{i-1}$.
By Proposition~\ref{pro5-1} \ref{pro5-1-no3} 
and \ref{pro5-1-no4}, $S_i=S_{i-1}\cup \{M-i\}$ 
is ideal in $D_{i}=(D_{i-1})_{u_i\rightarrow M-i}$.

Hence Claim 2 holds.

By Claim 2 and 
Proposition~\ref{pro5-0},
$\Re(D_i)=\Re(D_i-S_i)$ for all $i=0,1,\cdots,n-1$.
Hence $\Re(D_{n-1})=\Re(D[\{u_n\}])=\emptyset$.

Note that $u_j<M-(i+1)<M-i=\min S_i$ for all $j:i+1\le j\le n$.
Thus, by Claims 1, 2 and 
Proposition~\ref{pro5-2}, 
$D_{i}\succeq D_{i+1}$, 
as $D_{i+1}=(D_i)_{u_{i+1}\rightarrow M-(i+1)}$.
\proofend


\begin{cor}\relabel{th5-2-cor1}
 $\Re(D_{n-1})=\emptyset$ and 
\begin{equation}
\Psi(D,x)-\Psi(D_{n-1},x)=
\sum_{i=0}^{n-3} d_i {x+i\choose n-2},
\end{equation}
where $d_i\ge 0$ for all $i=0,1,\cdots,n-3$.
Furthermore, 
$\W(D)=\emptyset$ \iff $d_i=0$ for all 
$i=0,1,\cdots,n-3$. 
\end{cor}

\proof
By Theorem~\ref{theo5-2}, $\Re(D_{n-1})=\emptyset$.
By Theorem~\ref{theo5-2} again,
$D_i\succeq D_{i+1}$ for all $i=0,1,\cdots,n-2$.
By Proposition~\ref{prop4-0}, $D_0\succeq D_{n-1}$.
Thus, the result  follows 
from  Proposition~\ref{prop4-1}.
\proofend

{\it Proof of Theorem~\ref{th-new2}}:
As $\Re(D_{n-1})=\emptyset$, 
$\Psi(D_{n-1},x)=\Omega(\bar D_{n-1},x)$
by Proposition~\ref{stan1}.
As $\Omega(\bar D_{n-1},x)=\Omega(\bar D,x)$,
Theorem~\ref{th-new2} follows from 
Corollary~\ref{th5-2-cor1}.
\proofend
  
We end this section with a special sink-elimination 
ordering of $D$ determined by the injective mapping $L:V\rightarrow \{1,2,\cdots,n\}$
defined in Section 2.
Assume that $V=\{v_1,v_2,\cdots,v_n\}$
and $L(v_i)=n+1-i$ for $i=1,2,\cdots,n$.
Then $\alpha=(v_1,v_2,\cdots,v_n)$ is 
a sink-elimination ordering of $D$.
Assume that $\Gamma_{D,\alpha}=(D_0,D_1,\cdots,D_{n-1})$.
So $D_{n-1}$ is obtained from $D$ 
by relabeling each $v_i$ by $n+1-i+\max V$ 
for all $i=1,2,\cdots,n-1$.
Recall that $D_L$ denotes the digraph obtained 
from $D$ by relabeling each vertex $v_i$ by $L(v_i)=n+1-i$.
By definition, $D_{n-1}^*$ is exactly the digraph $D_L$,
implying that  $\Psi(D_L,x)=\Psi(D^*_{n-1},x)=\Psi(D_{n-1},x)$.
By Corollary~\ref{th5-2-cor1}, we have the following conclusion.

\begin{cor}
$\Psi(D,x)=\Psi(D_L,x)$ \iff 
$\W(D)=\emptyset$.
\end{cor}

\section{Proof of Theorem~\ref{cor-new1}
\relabel{sect-chro}
}

Let $G=(V,E)$ be a simple graph,
where $V=[n]$.
Recall that $\P(V)$ is the set of orderings  
of members of $V$.  
So $|\P(V)|=n!$.
Recall that 
$\AO(G)$ is the set of acyclic orientations of $G$.
Then $\P(V)$ can be partitioned according to members $D$ of  
$\AO(G)$ as stated in the following lemma.

\begin{lem}\relabel{le6-1} 
(i)  $\P(V)=\bigcup_{D\in \AO(G)}\OP(D)$;

(ii) $\OP(D_1)\cap \OP(D_2)=\emptyset$ for any pair of 
distinct orientations $D_1,D_2\in \AO(G)$.
\end{lem}

\proof (i). Clearly, $\OP(D)\subseteq \P(V)$.
For an ordering $\pi = (a_1, a_2,\cdots, a_n)$ of the elements of 
$V = [n]$,  
if $D$ is the orientation of $G$ such that 
$(a_i,a_j)\in A(D)$ whenever $i<j$ and $a_ia_j\in E$,
then $\pi\in \OP(D)$. 
Thus (i) holds.

(ii). Suppose that $D_1,D_2\in \AO(G)$ and 
$\pi=(a_1,a_2,\cdots,a_n)\in \OP(D_1)\cap \OP(D_2)$.
For any edge $a_ia_j$  in $G$, 
$i<j$ implies that $(a_i,a_j)\in A(D_1)\cap A(D_2)$.  
Thus $D_1$ and $D_2$ are the same.
Hence (ii) holds.
\proofend

\begin{lem}\relabel{le6-2} 
For any simple graph $G$, 
\begin{equation}\relabel{le6-2-eq1}
\Psi(G,x)=\sum_{D\in \AO(G)}\Psi(D,x).
\end{equation}
\end{lem}

\proof
Let $D\in \AO(G)$.
For any  $\pi=(a_1,a_2,\cdots,a_n)\in \OP(D)$
and any $i=1,2,\cdots,n-1$,
$a_i$ and $a_{i+1}$ are adjacent in $G$ \iff   
$a_i\rightarrow a_{i+1}$ in $D$, implying that 
$\delta_G(a_i,a_{i+1})=\delta_D(a_i,a_{i+1})$.
Thus $\delta_G(\pi)=\delta_D(\pi)$ holds
for any $\pi\in \OP(D)$, implying that 
\begin{equation}\relabel{cor-new1-eq2}
\Psi(D,x)=\sum_{\pi\in \OP(D)}{x+\delta_D(\pi)\choose n}
=\sum_{\pi\in \OP(D)}{x+\delta_G(\pi)\choose n}.
\end{equation}
Then, by Lemma~\ref{le6-1}, 
\begin{equation}\relabel{cor-new1-eq3}
\Psi(G,x)=\sum_{\pi\in \P(V)}{x+\delta_G(\pi)\choose n}
=\sum_{D\in \AO(G)}\sum_{\pi\in \OP(D)}{x+\delta_G(\pi)\choose n}.
\end{equation}
Thus (\ref{le6-2-eq1}) follows from (\ref{cor-new1-eq3}) 
and (\ref{cor-new1-eq2}).
\proofend

\begin{pro}\relabel{pro6-2} 
For any simple graph $G=(V,E)$ with $V=[n]$,
where $n\ge 3$,  
\begin{equation}\relabel{p6-2-e1}
\Psi(G,x)=(-1)^n \chi(G,-x)+
\sum_{i=0}^{n-3} d_i{x+i\choose n},
\end{equation}
where $d_i\ge 0$ for all $i=0,1,\cdots,n-3$.
Furthermore, $d_i=0$ for all $i=0,1,\cdots,n-3$
\iff $\W(G)=\emptyset$.
\end{pro}

\proof  By Theorem~\ref{stan2}, 
\begin{equation}\relabel{co-new1-eq1}
(-1)^n\chi(G,-x)=\sum_{D\in \AO(G)}\Omega(\bar D,x).
\end{equation}
Then, Lemma~\ref{le6-2}, (\ref{co-new1-eq1}) 
and Theorem~\ref{th-new2} imply that 
\begin{align}\relabel{pro6-2-eq2}
\Psi(G,x)-(-1)^n\chi(G,-x)
&=\sum_{D\in \AO(G)}
\big (\Psi(D,x)-\Omega(\bar D,x)\big )
\nonumber \\
&=\sum_{D\in \AO(G)}\sum_{i=0}^{n-3} 
d_{D,i}{x+i\choose n-2},
\end{align}
where $d_{D,i}\ge 0$ for all $i=0,1,\cdots,n-3$,
and $d_{D,i}=0$ for all $i=0,1,\cdots,n-3$
\iff $\W(D)=\emptyset$.
Thus $d_i=\sum\limits_{D\in \AO(G)}d_{D,i}\ge 0$
for all $i=0,1,\cdots,n-3$.
If $\W(G)=\emptyset$, then $\W(D)=\emptyset$ for all 
$D\in \AO(G)$, 
implying that $d_i=0$ for all $i=0,1,\cdots,n-3$.
Thus, it remains to show that 
if $\W(G)\ne \emptyset$, 
then $d_i>0$ for some $i$.

Now assume that  $\W(G)\ne \emptyset$.
Then  $\W(D)\ne \emptyset$ for some $D\in \AO(G)$.
By Theorem~\ref{th-new2}, $d_{D,i}>0$ for some $i$,
implying that $d_i>0$.

Hence Proposition~\ref{pro6-2} holds.
\proofend

Theorem~\ref{cor-new1} 
follows directly from Proposition~\ref{pro6-2}.

\section{Further study
\relabel{sec7}
}

We end this article with some problems 
that may merit further study.  
We assume that $G=(V,E)$ is a simple graph with $V=[n]$,
unless otherwise stated.

\subsection{Possible extensions of 
Theorems~\ref{cor-new1} and~\ref{th-new2}}

Theorem~\ref{cor-new1} gives an expression for 
$\chi(G,x)$ in terms of a summation over all $n!$  
orderings of elements of $V$
whenever $\W(G)=\emptyset$.
This result is established by applying Theorem~\ref{th-new2} 
and Stanley's result, Theorem~\ref{stan2}.
Is it possible to find new results analogous to  
Theorems~\ref{cor-new1} and~\ref{th-new2}
by revising $\delta_D(\pi)$ 
and $\delta_G(\pi)$ for  orderings $\pi$ of elements of $V$
such that the new results hold for a larger family of 
acyclic digraphs $D$ and a larger family of 
simple graphs $G$ respectively?

\subsection{Graphs $G$ with $\W(G)=\emptyset$
\relabel{sec7-1}}

Recall that 
$\W(G)=\{\{a,b,c\}: 1\le a<b<c\le n, ac\in E, ab,bc\notin E\}$.
By definition of $\W(G)$, the following observation follows
directly. Let $N_G(a)$ denote the set of vertices in $G$
which are adjacent to $a$.

\begin{pro}\relabel{pro7-1}
$\W(G)=\emptyset$ \iff 
for every edge $ac\in E$ with $a<c$,
$\{b: a<b<c\}\subseteq   N_G(a)\cup N_G(c)$ holds.
\end{pro}

For a bijection $\omega:V\rightarrow [n]$,
let $G_{\omega}$ be the graph obtained from $G$ 
by relabeling each vertex $v\in V$ by $\omega(v)$.
Let $\Gn$ denote the set of simple graphs $G=(V,E)$ 
such that $\W(G_{\omega})=\emptyset$ for some 
bijection $\omega: V\rightarrow [n]$.

If $G$ is a complete multi-partite graph and 
$ac$ is an edge in $G$, 
then $u\in N_G(a)\cup N_G(c)$ holds 
for every $u\in [n]-\{a,c\}$.
Thus, by Proposition~\ref{pro7-1}, 
$\W(G_{\omega})=\emptyset$ holds for an arbitrary 
bijection $\omega:V\rightarrow [n]$.
If $G$ is not a complete multi-partite graph,
this property does not hold.

The observations in the following proposition 
can be verified easily.

\begin{pro}\relabel{pro7-2} 
\begin{enumerate}
\item If $G=(V,E)$ is a complete multi-partite graph
and $\omega:[n]\rightarrow [n]$ is a bijection, 
then $\W(G_{\omega})=\emptyset$ holds and thus $G\in  \Gn$;
\item if $G$ is disconnected,    
then $G\in \Gn$ \iff 
each component of $G$ 
belongs to $\Gn$;
\item if $G\in \Gn$, then the subgraph of $G$ 
induced by any subset $S\subseteq V(G)$
belongs to $\Gn$.
\end{enumerate}
\end{pro}

By Proposition~\ref{pro7-1}, we have the following
 relation 
between $\W(G)$ and $\W(G-u)$, where $u\in V$.

\begin{pro}\relabel{pro7-3}
Let $u\in \{1,n\}$.
If $\{w: \min\{u,v\}<w<\max\{u,v\}\}\subseteq 
N_G(u)\cup N_G(v)$ holds for every $v\in N_G(u)$, then 
$\W(G)=\emptyset$ \iff $\W(G-u)=\emptyset$. 
\end{pro}

By Proposition~\ref{pro7-3}, the following corollary follows.

\begin{cor}\relabel{cor7-3}
Let $u\in \{1,n\}$.
If either 
$u=1$ and $N_{G}(u)=\{2,3,\cdots,k\}$,
or $u=n$ and $N_{G}(u)=\{k,k+1,\cdots,n-1\}$,
where $1\le k\le n$,
then $\W(G)=\emptyset$ \iff $\W(G-u)=\emptyset$. 
\end{cor}

Applying Proposition~\ref{pro7-3} 
or Corollary~\ref{cor7-3}, we find a family of 
graphs in $\Gn$ which are not complete multi-partite graphs.

\begin{pro}\relabel{pro7-4}
Let $G=(V,E)$ be a simple graph on $n$ vertices.
\begin{enumerate}
\item If $u$ is a vertex in $G$ such that $G-u$ 
is a complete multi-partite graph, then $G\in \Gn$;
\item Assume that $\{u,v\}$ is an independent set of $G$
such that $G-\{u,v\}$ is a complete multi-partite graph.
If either $\{u,v\}$ is a dominating set of $G$ 
or $N_G(u)\cap N_G(v)=\emptyset$, then $G\in \Gn$.
\end{enumerate}
\end{pro}

\proof (i). Assume that $k=|N_G(u)|$.
Let $\omega$ be a bijection from $V$ to $[n]$ 
such that $\omega(u)=1$ 
and $2\le \omega(w)\le k+1$ for all $w\in N_G(u)$.
As $G-u$ is a complete multi-partite graph, 
by Proposition~\ref{pro7-2} (i), 
$\W((G-u)_{\omega'})=\emptyset$,
where $\omega'$  is the mapping of $\omega$ 
restricted to $V-\{u\}$.
By Corollary~\ref{cor7-3},
$\W(G_{\omega})=\emptyset$
and so $G\in \Gn$.
Thus (i) holds.

(ii). We first the case that $\{u,v\}$ is a dominating 
set of $G$.
Assume that $d_G(u)=n_1$ and $d_G(v)=n_2$.
Then $|N_G(u)\cap N_G(v)|=n_1+n_2-n+2$,
$|N_G(u)\setminus N_G(v)|=n-2-n_2$ and $|N_G(v)\setminus N_G(u)|=n-2-n_1$.
Let $\omega$ be a bijection from $V$ to $[n]$ 
such that 

(a) $\omega(u)=1$, $\omega(v)=n$;

(b) $2\le \omega(w)\le n-1-n_2$ for all $w\in N_G(u)
\setminus N_G(v)$;

(c) $n-n_2\le \omega(w)\le n_1+1$ 
for all $w\in N_G(u)\cap N_G(v)$; and

(d) $n_1+2\le \omega(w)\le n-1$ 
for all $w\in N_G(v)\setminus N_G(u)$.

As $G-\{u,v\}$ is a complete multi-partite graph, 
by Proposition~\ref{pro7-2} (i), 
$\W((G-\{u,v\})_{\omega''})=\emptyset$,
where $\omega''$  is the mapping of $\omega$ 
restricted to $V-\{u,v\}$.
As $\omega$ satisfies the above conditions (a), (b), (c) and (d),
by Corollary~\ref{cor7-3},
$\W((G-u)_{\omega'})=\emptyset$
and $\W(G_{\omega})=\emptyset$.
Thus $G\in \Gn$.

Now consider the case that $N_G(u)\cap N_G(v)=\emptyset$.
Assume that $d_G(u)=n_1$ and $d_G(v)=n_2$.
Then $n_1+n_2\le n-2$.
Let $\omega$ be a bijection from $V$ to $[n]$ 
such that 

(a') $\omega(u)=1$, $\omega(v)=n$;

(b') $2\le \omega(w)\le n_1+1$ for all $w\in N_G(u)$; and 

(c') $n-n_2\le \omega(w)\le n-1$ 
for all $w\in N_G(v)$.

As $G-\{u,v\}$ is a complete multi-partite graph, 
by Proposition~\ref{pro7-2} (i), 
$\W((G-\{u,v\})_{\omega''})=\emptyset$,
where $\omega''$  is the mapping of $\omega$ 
restricted to $V-\{u,v\}$.
As $\omega$ satisfies the above conditions (a'), (b') and (c'),
by Corollary~\ref{cor7-3},
$\W((G-u)_{\omega'})=\emptyset$
and $\W(G_{\omega})=\emptyset$.
Thus $G\in \Gn$.
\proofend

In general, 
it seems not easy to determine all graphs in $\Gn$.
We now propose the following problem.

\begin{prob}\relabel{prob3} 
Characterize the family $\Gn$.
\end{prob}

As an example of studying Problem~\ref{prob3},
we now consider trees.
For any tree $T$ on $n$ vertices,
if $T$ is a star or a path, then it can be verified 
easily that $T\in \Gn$.

Now assume that $n\ge 5$. 
Let $T'$ denote the tree  obtained from $T$ by removing all
vertices of degree $1$.
If $T'$  is a path,
we can prove that $T\in \Gn$.
Assume that $T'$ is a path of order $k$: $u_1u_2\cdots u_k$.
Then $T$ is a tree obtained from $T'$ by adding 
$c_i$ new vertices and adding a new edge joining $u_i$ 
to each of them for all $i=1,2,\cdots,k$,
where $c_1,c_2,\cdots,c_k$ are some non-negative integers.
We first label each $u_i$ by $c_1+\cdots+c_i+i$.
Then we label the $c_i$ leaves adjacent to $u_i$ 
by numbers in the set 
$\{j: c_1+\cdots+c_{i-1}+i\le j\le c_1+\cdots+c_i+i-1\}$.
An example is shown in Figure~\ref{f5}.

\begin{figure}[htbp]
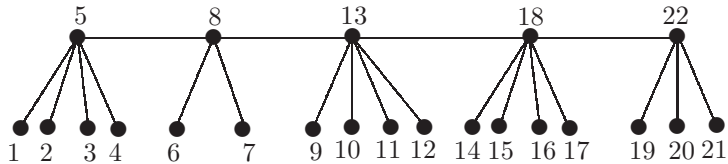

  \centering
\input f5.pic

  \caption{Labeling vertices in a tree $T$ whose non-leaf
  vertices induce a path}
\relabel{f5}
\end{figure}

If $T'$ is not a path, we believe $T\notin \Gn$.
For example, for the tree $T$ in Figure~\ref{f4},
$T'$ is not a path. It is left to the readers to 
verify that $T\notin \Gn$.

\begin{figure}[htbp]
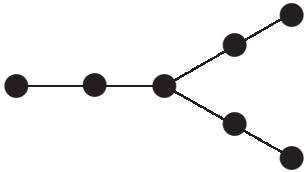

  \centering
\input f4.pic

  \caption{A tree which does not belong to $\Gn$}
\relabel{f4}
\end{figure}

\begin{con}\relabel{con7-1} 
For any tree $T$ on $n$ vertices, $T\in \Gn$ 
\iff either $n\le 2$ or the tree $T'$ obtained from $T$ 
by deleting all vertices of degree $1$ is a path.
\end{con}

\subsection{Interpretations of $d_i$'s in 
Theorem~\ref{th-new2}
\relabel{sec7-2}}

Let 
$D=(V,A)$ be an acyclic digraph with $V=[n]$, where $n\ge 3$, and 
let $d_0,d_1,\cdots,d_{n-3}$ be the numbers in 
Theorem~\ref{th-new2}.
In the following, we give an interpretation of $d_k$'s 
for a special case. 

\def \sink {{\cal S}ink}

Let $\sink(D)$ be the set of sinks of $D$.
Recall that for $u\in V$,
$R_D(u)$ is the set of $v\in V-\{u\}$ 
such that $D$ has a directed path from $u$ to $v$.
Thus  $u\in \sink(D)$ \iff $R_D(u)=\emptyset$.


\begin{pro}\relabel{pro5-5}
Let 
$u\in V$ with either $u\in \sink(D)$ 
or $\max (V\setminus R_D(u))<\min R_D(u)$ and 
$y<\min R_D(y)$ for each $y\in R_D(u)-\sink(D)$.
If $\W(D)=\W(D,u)$, 
then for $j=0,1,2,\cdots,n-3$, 
$d_j$ is the number of $D$-respecting orderings 
$\pi=(a_1,\cdots,a_i,u,a_{i+1},\cdots,a_{n-1})$
such that 

(i) either $\delta_D(\pi)=j+2$ and 
$\{a_i,u,a_{i+1}\}\in \W(D)$
with $a_i>a_{i+1}>u$; or

(ii) $\delta_D(\pi)=j$ and 
$\{a_i,u,a_{i+1}\}\in \W(D)$
with $a_i>u>a_{i+1}$.
\end{pro}

\proof Let $S=R_D(u)$. 
It is known that when $u\notin \sink(D)$, 
$\min S\ge 1+\max(V\setminus S)$.
Assume that $u\notin \sink(D)$ and $\min S= 1+\max(V-S)$.
Let $D'=(V',A')$ be the digraph obtained from $D$ 
by relabeling each vertex $v$  by $2v$
to get a new digraph $D'$. 
It is easy to verify that $W=\{a,b,c\}\in \W(D)$ 
\iff $W'=\{2a,2b,2c\}\in \W(D')$,
$\delta_{D}(\pi)=\delta_{D'}(\pi')$,
where $\pi'=(2a_1,2a_2,\cdots,2a_n)$ 
for $\pi=(a_1,a_2,\cdots,a_n)$,
and so $\Psi(D',x)=\Psi(D,x)$.
Observe that $\min S'\ge 2+\max(V'-S')$.
We can replace $D$ by $D'$ for the proof of this result.
 
Thus we may assume that either $u\in \sink(D)$
or $\min S\ge 2+\max(V\setminus S)$. 
We are now going to 
complete the proof by showing the following claims.

\inclaim Let $r=\max(V)+1$ if $u\in \sink(D)$, 
and let $r=\min (S)-1$ otherwise.
Then $\Psi(D_{u\rightarrow r},x)=\Omega(\bar D,x)$. 

Clearly, $r>\max \P_D(u)$ when $u\in \sink(D)$, 
and $\max\P_D(u)< r<\min F_D(u)$ otherwise.
It is also clear that $r\notin V$.
By Lemma~\ref{lem4-1}, $\W(D_{u\rightarrow r})=\W(D-u)$.
By the given condition, $\W(D-u)=\emptyset$,
implying that  $\W(D_{u\rightarrow r})=\emptyset$.
Thus, by Theorem~\ref{th-new1}, $\Psi(D_{u\rightarrow r}),x)
=\Omega(\bar D_{u\rightarrow r},x)=\Omega(\bar D,x)$.

\inclaim 
Let $r$ be the number given in the previous Claim.
Then 
$\Psi(D,x)-\Psi(D_{u\rightarrow r},x)
=\sum\limits_{j=0}^{n-3}d_j{x+j\choose n-2}$,
where $d_j$ is the number defined  
in Proposition~\ref{pro5-5}.

By Theorem~\ref{th4-2}, 
$\Psi(D,x)-\Psi(D_{u\rightarrow r},x)
=\sum\limits_{j=0}^{n-3}d_j{x+j\choose n-2}$  
holds with $d_j=c_{j+2}(D,u)+c'_{j}(D,u)$ for $j=0,1,\cdots, n-3$.
By definitions of $c_{j}(D,u)$ and $c'_{j}(D,u)$,
Claim~\thecountclaim\ holds.

By Claims 1 and 2, the result holds.
\proofend

Note that Proposition~\ref{pro5-5} gives 
an interpretation of $d_i$  
for a special case only.
In general, we would like to propose the following 
problem.

\begin{prob}\relabel{prob4}
Interpret the numbers $d_i$ in Theorem~\ref{th-new2}.
\end{prob}

\noindent {\bf Acknowledgment}: 
The author wishes to thank the referees 
for their very helpful comments and suggestions.

\end{document}